%% file: main.tex
\input{preamble}

\title{Border Rank Lower Bounds for Families of $\GL(V)$-invariant Tensors}
\author{Suhas Vadan Gondi}
\address{Department of Mathematics, University of California San Diego}
\email{sgondi@ucsd.edu}
\begin{document}

\begin{abstract}
    We give non-trivial lower bounds for the border rank of families of $\GL(V)$-invariant tensors in $U\otimes \bS_\lambda V\otimes \bS_\mu V$ where $U$ is $V$, $\Sym^2V$ or $\bigwedge^2V$. In particular, we provide a family of tensors with border rank reaching arbitrarily close to $2\ell$ in the unbalanced case, where $\ell$ is the largest ambient vector space dimension. We do this by resolving a conjecture introduced by Wu, and obtaining new results on $6j$-symbols as a byproduct. We then generalize our results to $\Sym^2V$ and $\bigwedge^2 V$ using novel techniques based on an application of a theorem of Kostant and Kempf collapsing.   
    
\end{abstract}
\maketitle
\section{Introduction}

\input{introduction}

\section{Preliminaries and Notation}\label{prelims section}
\input{prelims}

\section{$\GL(V)$-invariant tensors obtained from Pieri maps}\label{GL(V)Section}

\input{GLVinvariantTensors}

\section{Conditions on $\lambda$ and $\mu$ for $T_{\lambda\mu}'$ to be an isomorphism}\label{Full rank section}

\input{fullrankofmatrix}

\section{Generic Rank of Pieri Maps}\label{Pieri section}

\input{genericrank}

\section{Explicit Lower Bounds and Generalisations}\label{generalisation}
\input{generalisations}

\bibliographystyle{alpha}
\bibliography{reference}
\end{document}

%% file: preamble.tex
\documentclass[12pt]{amsart}

\usepackage{enumerate, amsmath, amsthm, amsfonts, amssymb, xy,  mathrsfs, graphicx, paralist, fancyvrb, ytableau, tikz-cd,epic,wasysym,upgreek}

\usepackage[margin=1in]{geometry} 
\usepackage[bookmarks, colorlinks=true, linkcolor=blue, citecolor=blue, urlcolor=blue]{hyperref}

\input xy
\xyoption{all}

\numberwithin{equation}{subsection}
\newtheorem{theorem}[equation]{Theorem}

\newtheorem{proposition}[equation]{Proposition}

\newtheorem{corollary}[equation]{Corollary}

\theoremstyle{definition}
\newtheorem{rmk}[equation]{Remark}
\newenvironment{remark}[1][]{\begin{rmk}[#1] \pushQED{\qed}}{\popQED \end{rmk}}
\newtheorem{eg}[equation]{Example}
\newenvironment{example}[1][]{\begin{eg}[#1] \pushQED{\qed}}{\popQED \end{eg}}
\newtheorem{defn}[equation]{Definition}
\newenvironment{definition}[1][]{\begin{defn}
[#1]\pushQED{\qed}}{\popQED \end{defn}}
\newtheorem*{mainresult}{Main Result}
\newenvironment{Main Result}[1][]{\begin{mainresult}[#1] \pushQED{\qed}}{\popQED \end{mainresult}}

\makeatletter
\@addtoreset{equation}{section}
\renewcommand{\thesubsection}{%
  \ifnum\c@subsection<1 \@arabic\c@section
  \else \thesection.\@arabic\c@subsection
  \fi
}
\makeatother

\newcommand{\bC}{\mathbf{C}}

\newcommand{\bL}{\mathbf{L}}

\newcommand{\cM}{\mathcal{M}}

\newcommand{\cO}{\mathcal{O}}

\newcommand{\cQ}{\mathcal{Q}}

\newcommand{\cR}{\mathcal{R}}

\newcommand{\bS}{\mathbf{S}}

\newcommand{\cU}{\mathcal{U}}
\newcommand{\fU}{\mathfrak{U}}

\newcommand{\cV}{\mathcal{V}}

\newcommand{\fg}{\mathfrak{g}}

\newcommand{\fp}{\mathfrak{p}}

\newcommand{\fu}{\mathfrak{u}}

\renewcommand{\phi}{\varphi}

\makeatletter
\def\Ddots{\mathinner{\mkern1mu\raise\p@
\vbox{\kern7\p@\hbox{.}}\mkern2mu
\raise4\p@\hbox{.}\mkern2mu\raise7\p@\hbox{.}\mkern1mu}}
\makeatother

\DeclareMathOperator{\Sym}{Sym}

\newcommand{\GL}{\mathbf{GL}}

\newcommand{\Gr}{\mathbf{Gr}}

\newcommand{\fsl}{\mathfrak{sl}}
\newcommand{\fgl}{\mathfrak{gl}}
\newcommand{\fso}{\mathfrak{so}}
\newcommand{\fsp}{\mathfrak{sp}}

%% file: introduction.tex
For complex vector spaces $A$, $B$ and $C$, the \textit{rank} of a tensor $T\in A\otimes B\otimes C$, denoted $\mathrm{Rk}(T)$, is defined to be the smallest $r$ such that $T$ can be written as the sum of $r$ elements of the form $a\otimes b\otimes c$. Computing the rank of an arbitrary tensor is NP hard \cite{HASTAD1990644}, but an important problem with applications in complexity theory, algebraic statistics and quantum information theory. For example, it is of practical importance to bound the rank of the matrix multiplication tensor $M_{uvw}\in \mathbb{C}^{uv}\otimes \mathbb{C}^{vw}\otimes \mathbb{C}^{wu},$ as this is directly linked to the complexity of the matrix multiplication operation, see for example \cite{Strassen1969, Conner_Harper_Landsberg_2023, v011a011}. 

The \textit{border rank} $\underline{\mathrm{Rk}}(T)$, which is easier to study than rank using algebraic techniques, is defined as the smallest $r$ such that $T$ is in the closure of the set of tensors of rank less than or equal to $r$. By definition, the set of tensors of border rank less than or equal to $r$ is an algebraic variety and it is an important but hard problem to obtain the defining equations of this variety. One reason why it is important to study these equations is to obtain lower bounds on the border ranks of tensors. One of the main results of our paper is the following.

\begin{mainresult}
    We construct families of tensors in $\mathbb{C}^\ell\otimes \mathbb{C}^m\otimes \mathbb{C}^n$, with $\ell>m>n$, having record-high border rank lower bounds. In particular, we exhibit a family whose border rank lower bound approaches $2\ell$. 
\end{mainresult}

Let $\ell,$ $m$ and $n$ denote the dimensions of $A$, $B$ and $C$ respectively. The expected generic border rank of a tensor in $A\otimes B\otimes C$, which has been proven for most values of $\ell$, $m$ and $n$ \cite[Section 3.1.4]{landsberg2012tensors},  is \begin{equation}\label{typicalRank} \left \lceil \frac{\ell mn}{\ell+m+n-2} \right \rceil.\end{equation} If the generic border rank is not the expected, then it is known to be higher. Thus, the border rank of a randomly chosen tensor in $A\otimes B\otimes C$ is at least \eqref{typicalRank} with probability one. However, barring some special cases, there are no known examples of explicit tensors whose border rank is proved to be the generic border rank and all known examples of lower bounds fall far short of the generic rank. Thus, this gives rise to an instance of a \textit{hay in a haystack} problem \cite{Hartnett2018}. 

The balanced case when $\ell=m=n$ holds is the situation that has the most literature. The generic rank in this case simplifies to about $\frac{\ell^2}{3}$. the best known family of tensors with high border rank for the balanced case are due to Landsberg and Michałek, who write down an explicit sequence of tensors in $\mathbb{C}^\ell\otimes\mathbb{C}^\ell\otimes\mathbb{C}^\ell$ for $\ell$ sufficiently large, having border rank at least $2.02\ell$ \cite{v021a013}. All previously known families of tensors in the balanced case did not breach the $2\ell$ bound for the lower bound on border rank. In the semi-balanced case when $n<m=\ell$ holds, the matrix multiplication tensor $M_{uvv}$ provides a family of examples of explicit tensors with border rank lower bound close to $2\ell$ \cite{v011a011}. In the unbalanced case when all the dimensions are different, there were almost no known examples of families of tensors with high border rank prior to Wu's work \cite{Wu24}. In this paper, we bridge this gap by building on the techniques from \cite{Wu24} to prove that a family of examples of explicit tensors in the unbalanced case have high border rank. 

One of the main ways to obtain lower bounds on border rank is to show that the tensor does not vanish on one of the defining equations of the corresponding secant variety to the Segre variety. This is the basic idea behind our approach, as we use novel techniques from representation theory and algebraic geometry to find the relevant defining equations and the secant varieties they are attached to and show that our tensors do not vanish on these equations. \textit{Flattenings} serve as a rich source of  equations, the idea behind which is to convert a tensor into a matrix whose rank can be obtained by computing its minors. Young flattenings are a class of flattenings introduced by Landsberg and Ottaviani \cite{LandsbergOttaviani2013}. They use a specific type of Young flattening, known as Koszul flattening, to improve the border rank lower bound of the matrix multiplication tensor \cite{v011a011}. Our approach uses Young flattenings to obtain lower bounds on the tensors we study. We show non-vanishing of these equations by proving invertibility of the flattening maps.      

The tensors we deal with originate from \cite[Section 4]{landsberg2022equivariantspacesmatricesconstant} and were studied by Wu in \cite{Wu24}. Let $V$ be a complex vector space of dimension $n$ and fix partitions $\lambda$ and $\mu$ such that $\mu$ is obtained from $\lambda$ by adding a box. We have the Pieri map (Section \ref{Pieri Rule}) $$\bS_{\mu}V\to \bS_{\lambda}V\otimes V, $$ which defines a $\GL(V)$-invariant tensor $T_{\lambda,\mu}\in V\otimes \bS_{\lambda}V\otimes \bS_{\mu}V^*$. Wu constructed the following Young flattening $T_{\lambda \mu}'$ for these tensors (Section \ref{GL(V)Section}):
  $$\bS_\lambda V\otimes \bS_\mu V\to \bS_\lambda V\otimes V\otimes \bS_{\lambda} V \to \bS_\mu V\otimes \bS_\lambda V,$$
  from which he obtained the following lower bound on border rank:
  $$\mathrm{\underline{Rk}}(T_{\lambda,\mu})\geq \left\lceil \frac{\mathrm{rk}(T_{\lambda,\mu}')}{r}\right\rceil ,$$ where $r$ is the rank of a non-zero matrix in the image of $V\to \mathrm{Hom}(\bS_{\lambda}V,\bS_{\mu}V)$, and we use $\mathrm{rk}(M)$ to denote the rank of a matrix. For the lower bound to be a good lower bound, we want $\mathrm{rk}(T_{\lambda,\mu}')$ to be large and $r$ to be small. It is shown in \cite[Section 2.2]{Wu24} or in \cite[Section 4.1]{landsberg2022equivariantspacesmatricesconstant} that $r$ is not of maximal rank if $\mu$ is obtained from $\lambda$ by adding the box to any row except the first or $n$th row. In this paper, we give formulas for $r$ that give explicit non-trivial lower bounds on border rank. As for the rank of $T_{\lambda \mu}'$, Wu conjectured \cite[Conjecture 1]{Wu24} that the Young flattening is an isomorphism for all pairs of partitions $\lambda$ and $\mu$. The first main result of this paper is to give a proof of a generalization of this conjecture in Section \ref{Full rank section}. 

\begin{theorem}\label{d=1Theorem}
 Suppose $\lambda$ and $\mu$ are two partitions such that $\mu$ is obtained from $\lambda$ by adding $d$ boxes with all boxes added to the same row (resp. same column). Then, the map obtained by composing the Pieri projection on the second factor with the Pieri inclusion on the first factor: $$ 
\bS_{\lambda}V\otimes \bS_{\mu}V\to \bS_{\lambda}V\otimes U\otimes \bS_{\lambda}V\to \bS_{\mu}V\otimes \bS_{\lambda}V,$$ 
is an isomorphism, where $U=\Sym^d V$ (resp. $\bigwedge^d V$). In particular, when $d=1$, it is an isomorphism for all partitions $\lambda$ and $\mu$ that differ by a box. 
\end{theorem}

In Section \ref{6jSection}, we define and show how this result is directly linked to $6j$-symbols and Racah W-coefficients, which are important topics in mathematical physics and representation theory. See \cite[Section 9]{venkatesh2026tetrahedralor6jsymbol} or \cite{Roberts:1998zka} for their role in mathematics and \cite{Martin_1960} or \cite{Varshalovich1989} for their role in physics. Theorem \ref{d=1Theorem} implies invertibility of the $6j$-symbols $$\begin{Bmatrix}
         \lambda & (d) & \mu \\ \lambda & \nu & \mu 
     \end{Bmatrix}\, , \begin{Bmatrix}
         \lambda & (1^d) & \mu \\ \lambda & \nu & \mu 
     \end{Bmatrix}$$ for all $\bS_\nu V$ in $\bS_\lambda V\otimes \bS_\mu V$.  

Motivated by the generality of Theorem \ref{d=1Theorem}, we next investigate the cases $U=\Sym^2 V$ or $\bigwedge^2 V$ when $\mu$ is obtained from $\lambda$ by adding two boxes to the same row or column depending on $U$. The Pieri map gives invariant tensors in $U\otimes \bS_\lambda V\otimes \bS_\mu V^*$, and we obtain similar invertible flattening maps for these tensors. The lower bound extracted in this case is the same as before with $r$ now denoting the rank of a generic matrix in the image of $$U\to \mathrm{Hom}(\bS_\lambda V,\bS_{\mu} V).$$ This is considerably harder than the previous case since the image is no longer a space of constant rank matrices. We use a theorem of Kostant to construct families of pairs of partitions $(\lambda, \mu)$ for which $r$ is not maximal and obtain non-trivial lower bounds on the border rank. The main idea is to use the resolutions spit out by Kostant's theorem to bound the rank of the cokernel of a certain map between modules over $\Sym U$, and prove that the rank of this cokernel equals $r.$ To the best of our knowledge, these techniques have not previously been used to study border rank of tensors. We hope that they will also prove useful for other problems related to tensor rank.

The generalizations to the $\Sym^2 V$ and $\bigwedge^2 V$ cases are interesting for several reasons. They extend the known families of tensors with high border rank (allowing for greater disparities among the dimensions of the ambient vector spaces than in the \(V\) case). They untilise new Young flattenings and help evaluate their stength in providing good lower bounds. Moreover, they address an interesting problem in their own right. Finally, these generalizations are interesting from another perspective, namely that of linear spaces of matrices of bounded rank. These are subspaces of matrices that contain no matrix of full rank, and they have been studied extensively \cite{HuangLandsberg2026, EISENBUD1988135, AtkinsonLloyd}. Our result provides conditions on the partitions $\lambda$ and $\mu$ that differ by two boxes such that the image $U\to \mathrm{Hom}(\bS_\lambda V,\bS_\mu V)$ is a linear space of matrices of bounded rank. 

We record the results obtained in the cases $U=V,$ $\Sym^2 V$ and $\bigwedge^2 V$ in the following theorem and refer the reader to Theorem \ref{ExplicitBoundTheorem} for the explicit lower bounds. We record the lower bounds for a selection of families of tensors, including the family with border rank lower bound close to $2\ell$ in Table \ref{tab:my_label}. 

\begin{theorem}\label{d=2Theorem}
Let $\lambda$ and $\mu$ be partitions. Then, we obtain non-trivial lower bounds for $\underline{\mathrm{Rk}}(T_{\lambda \mu})$ in the following situations.
\begin{enumerate}
    \item When $U=V,$ and $\mu$ is obtained from $\lambda$ by adding a box to any row except the first or $n$th row. 
    \item When $U=\Sym^2 V$, and $\mu$ is obtained from $\lambda$ by adding two boxes to the second row, and $\mu$ satisfies 
       $$\mu_1+\mu_3+3\leq 2\mu_2,$$
       and $\dim \bS_\lambda V\leq \dim \bS_\mu V$ holds. 

\item When $U=\bigwedge^2 V$, and $\mu$ is obtained from $\lambda$ by adding two boxes to the second and third rows but in the same column, and $\mu$ satisfies 
       $$\mu_1+\mu_4+2\leq 2\mu_2,$$
       and $\dim \bS_\lambda V\leq \dim \bS_\mu V$ holds. 
\end{enumerate}
\end{theorem}

\begin{corollary}
The border rank lower bounds for some families of $\GL(V)$-invariant tensors are listed in Table \ref{tab:my_label}. Here, $n$ denotes the dimension of $V$ and $\ell$ denotes the dimension of the largest ambient vector space. We note that the last column presents the limit of the lower bound obtained as $n\to \infty$. 

     \vspace{0.1in}
     \begin{center}
    \begin{tabular}{|c|c|c|c|c|} 
    \hline
         $U$&$\lambda$&$\mu$&Relation between $n$ and $d$ & Lower Bound for $\underline{\mathrm{Rk}}(T_{\lambda \mu})$  \\
         \hline 
         $V$& $(1^d)$& $(1^{d+1})$ & $n=2(d+1)$&$2\ell$\\
         \hline
          $V$& $(d,d-1)$& $(d,d)$ & $n=d+2$&$1.33\ell$\\
         \hline
         $\Sym^2 V$& $(d,d-2)$& $(d,d)$ & $n\sim (\sqrt 3-1) d$&$1.42\ell$\\
         \hline
          $\bigwedge^2 V$& $(d,d-1,d-1)$& $(d,d,d)$ & $n\sim (\sqrt 3-1) d$&$1.42\ell$\\
          \hline

    \end{tabular}\label{tab:my_label}
         \end{center}
    \vspace{0.1in}
 \end{corollary}

In Section \ref{generalisationSection}, we give a proof of the fact that the $2\ell$ lower bound can never be breached by our Young flattening in the $V$ case (Proposition \ref{2lCannotBreachProp}). Thus, we must look elsewhere for stronger bounds. Table \ref{tab:my_label} points to the fact that we get better bounds when looking at the $\Sym^2 V$ and $\bigwedge^2 V$ cases compared to the $V$ case when the boxes are added to the top three rows. Naturally, the next thing to ask is what we can say in the $\Sym^2 V$ and $\bigwedge^2 V$ cases when the boxes are added to the lower rows. Kostant's theorem unfortunately does not help with this case. Instead, we use Kempf collapsing to prove the following result in Section \ref{generalisationSection}, which sheds some light in this direction. 

\begin{theorem}\label{weymangeometricTheorem}
  Assume that $W$ is an $n-1$ dimensional complex vector space and that $\ell=\dim \bS_\mu W\geq \dim \bS_\lambda W$ holds. Let $k\geq \mu_1$, $\lambda^k=(k,\lambda_1,\dots\lambda_n)$ and $\mu^k=(k,\mu_1,\dots,\mu_n)$. Suppose that the border rank lower bound for $T_{\lambda \mu}$ (over $W$) obtained using the Young flattening is $\ell C$. Then, the border rank lower bound for the family of tensors $T_{\lambda^k \mu^k}$ (over $V$) tends to $\ell C$ (rapidly) as $k$ goes to infinity.   
\end{theorem}
We note that this is a general result that holds for $\Sym^d V$ or $\bigwedge^d V$, and hence in particular provides a way to construct an infinite family of examples of tensors with high border rank over $\Sym^d V$ and $\bigwedge^d V$ from a single example.

Finally, we discuss potential paths and obstacles to generalising the above methods to $\Sym^d V$ and $\bigwedge^d V$ for $d>2$. We apply Kostant's theorem on $E_6$ to obtain families of pairs of partitions $(\lambda,\mu)$ such that a generic matrix in the image of $\bigwedge^3 \mathbb{C}^6\to \mathrm{Hom}(\bS_{\lambda} V,\bS_\mu V)$ is not of maximal rank.

\subsection{Acknowledgments} I am grateful to my advisor Steven Sam for useful discussions, for suggesting ideas and for reading numerous drafts of this paper and providing valuable feedback. I would also like to thank J. M. Landsberg and Derek Wu for introducing the problems in this paper to me and for helpful conversations.

%% file: prelims.tex
\subsection{Partitions, Schur Functors and Pieri Maps}\label{Pieri Rule}

 A partition $\lambda$ of $n$ is a decreasing tuple of non-negative integers that sum up to $n$. For a partition $\lambda=(\lambda_1,\lambda_2\dots)$ we denote the length $l(\lambda)$ to be the largest $i$ for which $\lambda_i\neq 0$. When a partition has repeated entries, we use exponential notation. For example, the partition $(2,2,1,1,1)=(2^2,1^3)$. We often denote a partition using its Young diagram. For example:  $$(3,2,1)=\tiny\ydiagram{3,2,1}$$
 
 Throughout this paper, let $V$ be a complex vector space of dimension $n$. There is a correspondence between the finite-dimensional irreducible algebraic representations of $\GL(V)$ and weakly decreasing integer sequences $\lambda$ with at most $n$ rows. We denote the irreducible representation corresponding to $\lambda$ by the Schur functor $\bS_{\lambda}V$. While it is hard to describe the construction of a general Schur functor, there are two cases that are easy to describe: 
$$\bS_{(d)}V=\Sym^d{V}\qquad \qquad \qquad \bS_{(1^d)}V=\bigwedge^d(V)$$ 
For a construction of the general Schur functor, we refer the reader to \cite[Chapter 2]{Wey03}. Note that the notation for $\bS_{\lambda}V$ in \cite{Wey03} is $\bL_{\lambda^T}V$. $\bS_{\lambda}V$ is compatible with base change and hence we can construct $\bS_{\lambda}\cV$ for any vector bundle $\cV$. 

The one dimensional representation $\bS_{(1^n)}V$ will be denoted by $\det V$. We may twist any representation of $\GL(V)$ by a power of $\det V$ using the formula $$(\det V)^k\otimes \bS_{(\lambda_1,\lambda_2,\dots)}V=\bS_{(\lambda_1+k,\lambda_2+k,\dots)} V$$

The dual of an irreducible representation of $\GL(V)$ is also irreducible and is given by the formula $$\bS_{(\lambda_1,\lambda_2,\dots, \lambda_n)}V^*=\bS_{(-\lambda_n,\dots,-\lambda_2,-\lambda_1)}V$$
We denote $(-\lambda_n,\dots,-\lambda_2,-\lambda_1)$ by $-\lambda^{\mathrm{opp}}$.

The decomposition into irreducibles of the tensor product of two Schur functors is given by the Littlewood-Richardson rule \cite[Theorem 2.3.4]{Wey03}. An important special case of this that we will need is the Pieri rule, which gives a decomposition of the tensor product 
 $$\bS_{\lambda}V\otimes \Sym^d V=\bigoplus_{\lambda'}\bS_{\lambda '}V $$ where the direct sum is over all $\lambda '$ obtained from $\lambda$ by adding $d$ boxes with no two boxes added to the same column. We have an analogous rule for the exterior power $$\bS_{\lambda}V\otimes \bigwedge ^d V=\bigoplus_{\lambda'}\bS_{\lambda '}V $$
 where the direct sum is over all $\lambda '$ obtained from $\lambda$ by adding $d$ boxes with no two boxes added to the same row. From this it follows that there exist unique (up to scalar) $\GL$-invariant maps $$\bS_{\lambda}V\otimes \Sym^d V\to \bS_{\lambda '}V \qquad \qquad \bS_{\lambda '}V\to \bS_{\lambda}V\otimes \Sym^d V$$ known as \textbf{Pieri maps}. There are analogous maps for $\bigwedge^d V$.

 \begin{example} 
 We have
      $$\tiny\ydiagram{2,1}\otimes \tiny\ydiagram{2}=\tiny\ydiagram{4,1}\oplus \tiny\ydiagram{3,2}\oplus \tiny\ydiagram{3,1,1}\oplus \tiny\ydiagram{2,2,1} $$
      where we denote the Schur functor $\bS_\lambda V$ by the Young diagram of its corresponding partition.  
 \end{example}

 \subsection{Hook-Length Formula} Let $\lambda$ be a partition. The hook-length formula is a formula that computes the dimension of $\bS_\lambda V$. Before we state the formula, we need to define hook length.
 \begin{definition}
     For a partition $\lambda$, the hook $H_\lambda(i,j)$ of the box in the $i$th row and $j$th column of its Young diagram is defined to be the set of boxes in the Young diagram whose positions $(i',j')$ satisfy $i'=i$ and $j'\geq j$ or $i'\geq i$ and $j'=j$. The hook length is then defined to be $$h_\lambda(i,j)=|H_\lambda(i,j)|.$$ For example, the hook lengths of all the boxes for the partition $(4,2,1)$ are given below. 
     $$ \begin{ytableau}
         6 & 4 & 2 & 1\\
         3 & 1\\
         1
     \end{ytableau}$$
 \end{definition}
 \begin{proposition}\cite[Exercise 6.4]{FH13}
     $$\dim(\bS_\lambda V)=\prod_{i,j}\frac{n-i+j}{h_{\lambda}(i,j)}$$
     where the product is over all the pairs $(i,j)$ that number the row and column of the boxes of $\lambda$, and $h_\lambda(i,j)$ is the hook length of the corresponding box.
 \end{proposition}

\subsection{Borel-Weil-Bott} Let $X=\Gr(n-1,n)$ denote the Grassmanian consisting of hyperplanes inside the $n$ dimensional complex vector space $V$. We have the following tautological exact sequence of vector bundles over $X$ 
$$0\to \cR\to V\otimes \cO_X\to \cQ\to 0$$ where $\cR$ is the rank $n-1$ tautological subbundle and $\cQ$ is the tautological quotient bundle. We note that $\cQ$ is a line bundle and is isomorphic to $\cO(1)$. 

Let $S_n$ denote the symmetric group (the Weyl group of $\GL(V)$). Any permutation in $S_n$ can be written as a product of adjacent transpositions. For a permutation $\sigma\in S_n$, we define the length $l(\sigma)$ to be the minimum number of adjacent transpositions required to write $\sigma$ as a product of transpositions. It is known that $$l(\sigma)= \#\{i<j|\, \sigma(i)>\sigma(j)\}.$$  

Let $\lambda=(\lambda_1,\lambda_2,\dots)$ be a weakly decreasing sequence of integers and let $\mu=(d,\lambda_1,\lambda_2\dots)$ and $\rho=(n-1,n-2,\dots,0)$. For a permutation $\sigma\in S_n$, define the dotted action $$\sigma \bullet\mu=\sigma(\mu+\rho)-\rho.$$
\begin{theorem}(Borel-Weil-Bott, \cite[Corollary 4.1.9]{Wey03})\label{BWB} 

    One of the following two mutually exclusive cases occur:
    \begin{enumerate}
        \item There exists $\sigma\neq id$ such that $\sigma \bullet \mu=\mu$. Then, all the cohomology groups $H^*(X,\bS_\lambda \cR(d))$ vanish.
        \item There exists a unique $\sigma\in S_n$ such that $\sigma \bullet \mu=\nu$ is weakly decreasing. Then we have $$H^{l(\sigma)}(X,\bS_{\lambda}\cR(d))=\bS_{\nu}V$$ and all other cohomology groups vanish.

    \end{enumerate}
\end{theorem}

\subsection{Kostant's Theorem} We state Kostant's theorem following \cite{Kum}. We refer the reader to \cite[Chapter 1]{Kum} for a review of the definitions and concepts involved. Since we will mostly work with finite-dimensional semisimple Lie algebras in this paper (specifically, the families  $C_n$ and $D_n$), we state the theorem for this case even though it holds in general for any symmetrizable Kac-Moody algebra. 

Let $\fg$ be a finite-dimensional semisimple Lie algebra. Let $\fp_Y$ be the standard parabolic subalgebra corresponding to a subset $Y$ of simple roots. Let $\fg_Y$ denote the Levi subalgebra and $\fu_Y$ the nilpotent radical of $\fp_Y$. Let $L(\lambda)$ and $L_Y(\mu)$ denote the finite-dimensional highest-weight $\fg$-module and $\fg_Y$-module respectively, where $\lambda$ is a dominant weight for $\fg$ and $\mu$ is a dominant weight for $\fg_Y$. 

Let $W$ be the Weyl group of $\fg$ and $W_Y$ the Weyl group of $\fg_Y$. We have the length function $l(\sigma)$ on $W$ defined to be the minimum number of simple reflections required to write $\sigma$ as a product of simple reflections. Let $W_Y'$ denote the set of unique minimal length representatives in each coset of $W/W_Y$. If $\rho$ denotes the sum of fundamental weights of $\fg$, we have the dotted action on the set of weights as defined in the previous section. We note that $\sigma\in W_Y'$ if and only if $\sigma^{-1}\bullet \lambda$ is a dominant weight for $\fg_Y$ whenever $\lambda$ is a dominant weight for $\fg$. 

\begin{theorem}(Kostant, \cite[Theorem 3.2.7]{Kum}\label{Kostant})
    
Let $\lambda$ be a dominant weight of $\fg$. We have the following isomorphism as $\fg_Y$-modules 
$$H_p(\fu_Y,L(\lambda)^*)\cong \bigoplus_{\substack{w\in W_Y'\\l(w)=p}}L_Y(w^{-1}\bullet \lambda)^*.$$

\end{theorem}

We remark that in \cite[Theorem 3.2.7]{Kum}, the theorem is stated for $\fu_Y^-$, which is isomorphic to the dual of $\fu_Y$ as a representation of $\fg_Y$ and that is why we take duals in our statement. 

\begin{remark}\label{KostantRemark}
    Lie algebra homology has more than one equivalent definition, but the most useful to us will be the following definition $$H_n(\fg, M)=Tor_n^{\fU_\fg}(\bC,M),$$ where $\fU_\fg$ is the universal enveloping algebra of $\fg$ and $M$ is a $\fU_\fg$-module. We note that there is a correspondence between $\fg$-representations and $\fU_\fg$-modules. When $\fg$ is abelian, $\fU_\fg$ is the polynomial ring over $\fg$ and $H_n(\fg,M)$ is the $nth$ syzygy module of $M$. We use this fact to construct free resolutions.
\end{remark}
 \subsection{Border Rank of Tensors} Let $A$, $B$, $C$ be complex vector spaces of dimensions $\ell$, $m$ and $n$ respectively. A tensor $T \in A\otimes B \otimes C$ is said to be pure if it is of the form $a\otimes b\otimes c$.
\begin{definition} 
    The rank of a tensor $T\in A\otimes B\otimes C$ (denoted $\mathrm{Rk}(T)$) is the smallest $r$ such that $T$ can be written as the sum of $r$ pure tensors. 
\end{definition}
The set of tensors of rank at most $r$ is not closed in the Zariski topology (\cite[Section 2.4.5]{landsberg2012tensors}). This motivates the definition of the border rank of a tensor.

\begin{definition}
    The border rank $\mathrm{\underline{Rk}}(T)$ is the smallest integer $r$ such that $T$ is the limit of tensors of ranks less than or equal to $r$.
\end{definition}
It follows trivially that $$\mathrm{\underline{Rk}}(T)\leq \mathrm{Rk}(T)$$

Matrices can be thought of as two-tensors, and it can be shown that the rank (and border rank) as tensors equals the matrix rank. While it is hard to compute the rank of a general tensor, we know how to compute the rank of a matrix. The idea behind flattenings is to ``flatten" our tensor into a matrix and use the matrix rank to obtain bounds on the tensor rank. 

Given a tensor $T\in A\otimes B\otimes C$, we obtain three natural flattenings of $T$ into matrices: $$T_A\colon A^*\to B\otimes C,$$
$$T_B\colon  B^*\to A\otimes C,$$
$$T_C\colon C^*\to A\otimes B.$$
Since rank is subadditive, we obtain lower bounds for the border rank of $T$
$$\max\{\mathrm{rk}(T_A),\mathrm{rk}(T_B), \mathrm{rk}(T_C)\}\leq \mathrm{\underline{Rk}}(T).$$
A tensor $T$ is said to be \textbf{concise} if $T_A$, $T_B$ and $T_C$ are full-rank. If a tensor is not concise, we may view it as a tensor in the subspace of one of the vector spaces where it is now concise. Thus, it suffices to study concise tensors. The previous lower bound for concise tensors gives  $$\max\{\ell,m, n\}\leq \mathrm{\underline{Rk}}(T).$$ If equality holds above, we say that $T$ is of minimal border rank.

 \subsection{Young Flattenings}\label{YFsection}
 Fix partitions $\lambda$, $\lambda '$, $\mu$, $\mu'$, $\nu$, $\nu '$ where $\lambda'$, $\mu'$ and $\nu'$ are obtained from $\lambda$, $\mu$ and $\nu$ respectively by adding a box.

Define the Young flattening $T'$ to be the matrix obtained as the composition: 

$$\bS_{\lambda}A\otimes \bS_{\mu} B\otimes \bS_{\nu'}C^*$$
$$\downarrow $$
$$\bS_{\lambda}A\otimes \bS_{\mu}B\otimes \bS_{\nu}C^*\otimes C^*$$
$$\quad \, \, \, \, \, \, \, \, \downarrow (T_C)$$
$$\bS_{\lambda}A\otimes \bS_{\mu}B\otimes \bS_{\nu}C^*\otimes A\otimes B$$ 
$$\downarrow$$
$$\bS_{\lambda'}A\otimes \bS_{\mu'}B\otimes \bS_{\nu}C^*$$

Let $r_1$, $r_2$, $r_3$ be the ranks of a generic matrix in the images $$A\to \mathrm{Hom}(\bS_{\lambda}A,\bS_{\lambda'}A),$$ $$B\to \mathrm{Hom}(\bS_{\mu}B,\bS_{\mu'}B),$$ $$C^*\to \mathrm{Hom}(\bS_{\nu}C^*,\bS_{\nu'}C^*).$$

If $T$ has rank one, then the rank of $T'$ is bounded above by $r_1r_2r_3$. Thus, by subadditivty of rank, we get a lower bound for the border rank of $T$: 

\begin{equation}\label{BorderLowerBound}
  \mathrm{\underline{Rk}}(T)\geq \left\lceil \frac{rk(T')}{r_1r_2r_3}\right \rceil.  
\end{equation}

%% file: GLVinvariantTensors.tex
Fix partitions $\lambda$ and $\mu$ where $\mu$ is obtained from $\lambda$ by adding a box. We have the $\GL(V)$-equivariant map  $$f\colon \bS_{\mu}V\to \bS_{\lambda}V\otimes V,$$ which corresponds to a tensor $$T_{\lambda,\mu}\in V\otimes \bS_{\lambda}V\otimes \bS_{\mu}V^*.$$ This paper focuses on giving non-trivial lower bounds for the border rank of such tensors.

Using the notation from section \ref{YFsection}, we consider the Young flattening with parameters $\lambda=\lambda$, $\lambda'=\mu$, $\mu=\varnothing$, $\mu'=(1)$, $\nu=\varnothing$, $\nu'=(1)$ to obtain $T_{\lambda,\mu}'$:
    $$\bS_\lambda V\otimes \bS_\mu V\xrightarrow{1\otimes f}\bS_\lambda V\otimes V\otimes \bS_{\lambda} V\xrightarrow{g\otimes 1}\bS_\mu V\otimes \bS_\lambda V$$
where $g$ is the Pieri map in the other direction $$g\colon \bS_\lambda V\otimes V\to \bS_\mu V.$$ 
\eqref{BorderLowerBound} gives us the following lower bound on the border rank of $T_{\lambda,\mu}$: $$\mathrm{\underline{Rk}}(T_{\lambda,\mu})\geq \left\lceil \frac{\mathrm{rk}(T_{\lambda,\mu}')}{r}\right\rceil $$ where $r$ is the rank of a generic matrix in the image of $V\to \mathrm{Hom}(\bS_{\lambda}V,\bS_{\mu}V).$ 

More generally, if $\lambda$ and $\mu$ are two partitions such that $\mu$ is obtained from $\lambda$ by adding $d$ boxes with at most one box added per row (column), then the Pieri rule (Skew Pieri rule) gives us a tensor $T_{\lambda,\mu}\in \Sym^dV\otimes \bS_{\lambda}V\otimes \bS_{\mu}V^*$ (respectively $T_{\lambda,\mu}\in \wedge^dV\otimes \bS_{\lambda}V\otimes \bS_{\mu}V^*$). We may consider the corresponding flattening maps $T_{\lambda,\mu}'$ 
\begin{equation}\label{FlatteningMap}
\bS_{\lambda}V\otimes \bS_{\mu}V\to \bS_{\lambda}V\otimes \Sym^d V\otimes \bS_{\lambda}V\to \bS_{\mu}V\otimes \bS_{\lambda}V,
\end{equation}

\begin{equation}\label{FlatteningMap'}
\bS_{\lambda}V\otimes \bS_{\mu}V\to \bS_{\lambda}V\otimes \wedge^d V\otimes \bS_{\lambda}V\to \bS_{\mu}V\otimes \bS_{\lambda}V,
\end{equation}
from which we can extract similar lower bounds for the tensor ranks as for the $V$ case. We note that $r$ in this case is the generic rank in the image of $\Sym^d V$ (respectively $\bigwedge^d V$) in $\mathrm{Hom}(\bS_{\lambda}V,\bS_{\mu}V)$. 

We get non-trivial lower bounds if the following two conditions hold: 
\begin{enumerate}
\item $T_{\lambda,\mu}'$ is full rank and hence an isomorphism. 
\item The generic rank $r$ is less than the dimensions of the ambient vector spaces.
\end{enumerate}
Suppose $\dim(\bS_{\lambda}V)=k$ and $\dim(\bS_{\mu}V)=\ell$. If both conditions above are satisfied, we get 
\begin{equation}\label{k-lLowerBound}
\mathrm{\underline{Rk}}(T_{\lambda,\mu})\geq \left \lceil\frac{k\ell}{r}\right \rceil>\max\{k,\ell\} \end{equation} Thus, the bound we get would be a non-trivial bound. The next sections will focus on giving conditions on $\lambda$ and $\mu$ for which the two statements hold.

%% file: fullrankofmatrix.tex
\subsection{Statement and Proof} In this section, we give sufficient conditions on $\lambda$ and $\mu$ for which the Young flattenings $$T_{\lambda \mu}'\colon \bS_{\lambda}V\otimes \bS_{\mu} V\to \bS_{\mu}V\otimes \bS_{\lambda}V$$ coming from \eqref{FlatteningMap} and \eqref{FlatteningMap'} are full rank. The main result will be the following theorem.

\begin{theorem}\label{FullrankTheorem}

    Suppose $\lambda$ and $\mu$ are two partitions such that $\mu$ is obtained from $\lambda$ by adding $d$ boxes with all boxes added to the same row or column. Then \eqref{FlatteningMap} and \eqref{FlatteningMap'} are isomorphisms. In particular, when $d=1$, \eqref{FlatteningMap} is an isomorphism for all partitions $\lambda$ and $\mu$ that differ by a box.  
\end{theorem}

Let us fix some notation. Unless stated otherwise, for the rest of this section $\lambda$ and $\mu$ will denote partitions such that $\mu$ is obtained from $\lambda$ by adding $d$ boxes to the same row or column. Let $U$ denote either $\Sym^dV$ or $\bigwedge^d V$ depending on whether the $d$ boxes are added to the same row or column.  

The following theorem from \cite[Theorem 4.2]{FLS} will be a key ingredient in the proof of Theorem \ref{FullrankTheorem}. 

\begin{theorem}[Ford, Levinson, Sam]\label{FLSThm}
Let $W$ be another complex vector space of dimension $n$ and let $R=\Sym(V\otimes W^*)$. The following unique (up to scalar) $\GL(V)\times \GL(W)$-equivariant $R$-linear map is injective:

    $$\bS_\lambda W \otimes \bS_\mu V \otimes R \to \bS_\mu W  \otimes \bS_\lambda V  \otimes R(d).$$

\end{theorem}
We note that in \cite[Theorem 4.2]{FLS}, the statement assumes that $\lambda$ and $\mu$ differ by only one box, which is the $d=1$ case. But in \cite[Remark 4.8]{FLS}, the authors mention that the proof can be generalized to arbitrary $d$. We elaborate on this further in the next remark. 

\begin{remark}
    For a partition $\alpha$, we say $(i,j)$ is an outer border square if $(i,j)\notin \alpha$ and either $(i-1,j-1)\in \alpha$ or $i=1$ or $j=1$ hold. 
    We say that two partitions $\alpha$ and $\beta$ differ by a connected border strip if $\beta$ can be obtained from $\alpha$ by adding a connected subset of outer border squares. Adding $d$ boxes to the same row or column are examples of this.

\begin{figure}[h]
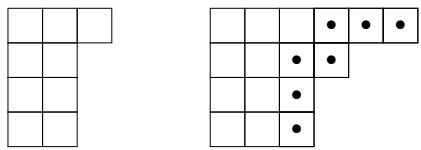

  \begin{center} 

\tiny \ytableausetup{centertableaux}
\ydiagram{3,2,2,2} \qquad \qquad 
\ytableausetup{nobaseline}
\ytableaushort
{\none \none \none \bullet \bullet \bullet,\none \none \bullet \bullet, \none \none \bullet, \none \none \bullet}
* {6,4,3,3}
\caption{$\alpha=(3,2^3)$ and $\beta=(6,4,3^2)$ are examples of partitions that differ by a connected border strip.}  
\label{figure1}
\end{center}
\end{figure}

According to \cite[Remark 4.8]{FLS}, Theorem \ref{FLSThm} still holds when $\lambda$ and $\mu$ differ by a connected border strip. However, a major obstruction in the general case is that it is possible that there are multiple $\GL(V)\times \GL(W)$-equivarant maps between the two spaces and Theorem \ref{FLSThm} gives injectivity for one of the maps which may or may not be the map we are interested in. Luckily, when $\lambda$ and $\mu$ differ by a horizontal or vertical strip there is a unique (up to scalar)  $\GL(V)\times \GL(W)$-equivarant map and we don't run into any problems. To see this, we note that it suffices to show that there is a unique $\GL(V)\times \GL(W)$-equivariant linear map between the degree $0$ components of the $R$-modules. Suppose the partitions differ by a horizontal strip (the same proof works for a vertical strip). Then, we observe that there is a unique chain of partitions $$\lambda=\lambda_0\subset\lambda_1\dots \subset \lambda_d=\mu,$$ where $\lambda_i$ differs from $\lambda_{i-1}$ by a box. Pieri rule says that there is a unique (up to scalar) equivariant map $\bS_{\lambda_i} V\to \bS_{\lambda_{i-1}}V\otimes V$ and taking compositions, there is a unique $\GL(V)$-equivariant map $$\bS_{\mu}V\to \bS_\lambda V\otimes V^{\otimes d}.$$
Since the Pieri map $$\bS_\mu V\to \bS_\lambda V\otimes \Sym^d V\to \bS_\mu V\otimes V^{\otimes d} $$ gives one such map, the above map has to be this. Similarly, there is a unique $\GL(W)$-equivariant map that factors as $$\bS_\lambda W\to \bS_\mu W\otimes \Sym^d W^*\to \bS_\mu W\otimes (W^*)^{\otimes d}$$ and combining the two, we observe that there is a unique $\GL(V)\times \GL(W)$-equivariant map that factors as 
\begin{equation}\label{factormap}
    \bS_\mu V\otimes \bS_\lambda W\to \bS_\lambda V\otimes \bS_\mu W\otimes \Sym^d V\otimes \Sym^d W^*\to \bS_\lambda V\otimes \bS_\mu W\otimes (V\otimes W)^{\otimes d}
\end{equation}
as required. 
\end{remark}

We will not give the proof of Theorem \ref{FLSThm} here but refer the reader to \cite[4.1.4]{FLS} for a proof. However, we will give a brief sketch of the proof by illustrating an example.

\begin{example}\label{FullrankExample}
    Suppose $n=4$, $\lambda=(6,2^3)$ and $\mu=(6,4,2^2)$. In \cite{Eisenbud_2011}, the authors construct pure free resolutions over $S=\Sym(V)$ involving Schur functors as modules and Pieri maps as differentials. The key fact that will be useful to us is that if $\lambda$ and $\mu$ differ by a connected border strip, then there exists a resolution that contains $\bS_\lambda V\otimes S$ and $\bS_\mu V\otimes S$ as modules in the resolution (this can be done by completing the border strip and then using the construction in \cite[Section 3]{Eisenbud_2011}). In particular, for our case, completing the border strip gives us the partitions from Figure \ref{figure1} and the resolution looks like (after twisting) $$0\to \ytableausetup{smalltableaux}\tiny\ydiagram{6,4,3,3}(-2)\to \tiny\ydiagram{6,4,3,2}(-1)\to\tiny\ydiagram{6,4,2,2}\xrightarrow{*} \tiny\ydiagram{6,2,2,2}(2)\to \tiny\ydiagram{3,2,2,2}(5).$$
    Here, we denote $\bS_\nu V\otimes S(d)$ by $\nu(d)$.
    
    Let $X$ denote $\mathbb{P}(W^*)\cong \Gr(3,4)$. Recall that we have the tautological exact sequence over $X$. We sheafify the above complex by replacing $V$ with the rank-$n$ bundle $V\otimes \cO_X(-1)$ on $X$, base-change along the flat extension $\Sym(V\otimes \cO_X(-1))\to R$, and tensor throughout with $\bS_{(6,1,1)}\cR$ to get terms of the form $$\nu(d)=\bS_\nu(V)\otimes\cO_X(d)\otimes \bS_{(6,1,1)}\cR\otimes R.$$  Let $\cM$ denote the sheaf resolved by this acyclic complex. Applying Theorem \ref{BWB}, we observe that the cohomologies of all the terms except $\lambda$ and $\mu$ vanish in all degrees, and we have $$H^1(X, \mu)=\bS_\mu V\otimes \bS_{(5,1,1,1)}W\otimes R, \qquad H^1(X,\lambda(2))=\bS_\lambda V\otimes \bS_{(5,3,1,1)}W \otimes R.$$
    Thus, the $E_1$ page of the hypercohomology spectral sequence for this complex only has two non-zero terms and looks like 
    
$$\begin{tikzcd}
0 & 0 & 0                                                       & 0                                                             & 0 & H^0(\cM) \\
0 & 0 & {\bS_\mu V\otimes \bS_{(5,1,1,1)} W\otimes R} \arrow[r] & {\bS_\lambda V\otimes \bS_{(5,3,1,1)} W\otimes R} \arrow[rru] & 0 & H^1(\cM)
\end{tikzcd}$$
We note that the diagonal map is in the $E_2$ page. 

Thus, the map in the $E_1$ page is injective. Finally, twisting both modules by $\det(W)$, we see that the map $$ \bS_\mu V\otimes \bS_{\lambda} W\otimes R\to \bS_\lambda V\otimes \bS_{\mu} W\otimes R $$ is injective.
\end{example}

\begin{remark}
    In the example preceding this remark, if we tensor the complex with $\bS_{(6,4,1)}\cR $ instead of $\bS_{(6,1,1)}\cR$, we get non-zero cohomology only for the third and fourth terms from the right. Consequently, we get injectivity in the case when $\lambda=(6,4,2,2)$ and $\mu=(6,4,3,2)$ (they differ by a box). 
    
    Tensoring with $\bS_{(6,4,2)}\cR$ gives us injectivity for when $\lambda=(6,4,2^2)$ and $\mu=(6,2,3^2)$ (they differ by two boxes in the same column). We note that in this case, the required map appears in the $E_2$ page of the hypercohomology spectral sequence. 
\end{remark}

\subsection*{Proof of Theorem \ref{FullrankTheorem}.} 
 Let the $R$-module $M$ be the cokernel of the map in Theorem \ref{FLSThm}. We have  $$0\to \bS_\lambda W \otimes \bS_\mu V \otimes R \to \bS_\mu W \otimes \bS_\lambda V \otimes R(d)\to M\to 0.$$ $M$ has projective dimension one and hence the codimension of $M$ is bounded above by one. We remark that $M$ cannot be free since its rank equals zero. Thus, the support of $M$ is a hypersurface in the space of matrices in $V^*\otimes W$ that is $\GL(V)\times \GL(W)$-invariant and there is exactly one such hypersurface: the space of non-invertible matrices. Hence, after replacing $W$ with $V$ and specializing to the identity matrix, we see that the resulting map is surjective. Restricting the map to the degree zero part, we obtain a map $$\bS_\lambda V \otimes \bS_\mu V  \to \bS_\mu V \otimes \bS_\lambda V$$ that is surjective and hence of full rank. We now show that this is the required map. 

  From \eqref{factormap}, we know that the degree zero part of the map in Theorem \ref{FLSThm} factors as follows: 
 
\begin{center}
\begin{tikzcd}
\bS_{\lambda}V\otimes \bS_{\mu}V \arrow[rd, "f"] \arrow[r ] & \bS_{\mu}V\otimes \bS_{\lambda}V\otimes \Sym^d(V\otimes V^*)                 \\
                                                      & \bS_{\mu}V\otimes \bS_{\lambda}V\otimes U \otimes U^* \arrow[u, "i"]
\end{tikzcd}
 \end{center} 
where the map $i$ is the inclusion map and map $f$ is the tensor product of the Pieri maps

  $$f_1 \colon \bS_\mu V \to \bS_\lambda V \otimes U,$$ $$f_2 \colon \bS_\lambda V \to \bS_\mu V\otimes U^*.$$ 
  Let us fix an element $x\otimes y \in \bS_\mu V\otimes \bS_\lambda V$ and a basis $\{e_i\}_{1\leq i\leq n}$ for $V$ with dual basis $\{e_i^*\}_{1\leq i\leq n}$ for $V^*$. This gives a basis $\{e^\alpha\}$ for $U$ where $\alpha$ is an $n$-tuple of non-negative integers (when $U$ is $\bigwedge^dV$, each $\alpha_i\leq 1$) and $e^\alpha=e_1^{\alpha_1}e_2^{\alpha_2}\dots e_n^{\alpha_n}$. Similarly, we obtain the dual basis $\{(e^\alpha)^*\}$ for $U^*$.   
  Let $$f_1(x)=\sum_{\alpha} x_\alpha\otimes e^\alpha,$$ $$f_2(y)= \sum_{\beta}y_\beta\otimes (e^\beta)^*.$$ Then  we have $$f(x\otimes y)=\sum_{\alpha,\beta} x_\alpha\otimes e^\alpha\otimes y_\beta\otimes (e^\beta)^*.  $$ 
  
Now, restricting to the identity matrix is the same as mapping $e^\alpha\otimes (e^\beta)^*$ to $\delta_{\alpha \beta}$. Thus the image after restriction would be $$\sum_{\alpha} x_\alpha\otimes y_\alpha,$$ which is exactly the image of $x\otimes y$ under the maps \eqref{FlatteningMap} and \eqref{FlatteningMap'}. This completes the proof. \qed

\subsection{$6j$-symbols}\label{6jSection}
 The maps \eqref{FlatteningMap} and \eqref{FlatteningMap'} are closely related to an important concept in mathematical physics and representation theory known as $6j$-symbols, whose general definition can be found in \cite[Section 1.42]{etingof2015tensor}. Here, we define them for representations of $\GL(V)$ for convenience. Let $\bS_\alpha V$, $\bS_\beta V$ and $\bS_\gamma V$ be three irreducible representations of $\GL(V)$. Due to associativity of tensor product, we have an isomorphism \begin{equation}\label{associativityIsomorphism} \bS_\alpha V\otimes (\bS_\beta V\otimes \bS_\gamma V)\to (\bS_\alpha V\otimes \bS_\beta V)\otimes \bS_\gamma V.\end{equation} Let $H_{\rho \sigma}^{\nu}$ denote the multiplicity spaces $\mathrm{Hom}(\bS_\nu V, \bS_\rho V\otimes \bS_\sigma V).$ Then, we have $$\bS_\rho V\otimes \bS_\sigma V=\bigoplus_\nu H_{\rho \sigma}^\nu \otimes \bS_\nu V.$$  Further, we have $$\bS_\alpha V\otimes (\bS_\beta V\otimes \bS_\gamma V)=\bigoplus_{\rho, \sigma }H_{\beta \gamma}^\rho\otimes H_{\alpha \rho }^{\nu}\otimes\bS_\nu V,$$ $$(\bS_\alpha V\otimes \bS_\beta V)\otimes \bS_\gamma V=\bigoplus_{\sigma, \nu }H_{\alpha \beta}^\sigma\otimes H_{\gamma \sigma }^{\nu}\otimes\bS_\nu V. $$ For a fixed irreducible representation $\bS_\nu V$ in $\bS_\alpha V\otimes \bS_\beta V\otimes \bS_\gamma V,$ the isomorphism \eqref{associativityIsomorphism} induces an isomorphism on the $\bS_\nu V$-isotypic component: \begin{equation}\label{6jDefinition}\bigoplus_{\rho}H_{\beta \gamma}^\rho\otimes H_{\alpha \rho }^{\nu}\to \bigoplus_\sigma H_{\alpha \beta}^\sigma\otimes H_{\gamma \sigma }^{\nu}.\end{equation} For six fixed parameters (and hence the six in the name), the $6j$-symbol $$\begin{Bmatrix}
     \alpha & \beta & \sigma \\ \gamma & \nu & \rho     
 \end{Bmatrix}$$ is defined as the component map in \eqref{6jDefinition}. 

 It is easy to see using the description of the maps \eqref{FlatteningMap} and \eqref{FlatteningMap'} that these maps are an isomorphism if and only if each of the $6j$-symbols $$\begin{Bmatrix}
     \lambda & (d)& \mu \\ \lambda & \nu & \mu
 \end{Bmatrix}\, , \begin{Bmatrix}
     \lambda & (1^d)& \mu \\ \lambda & \nu & \mu
 \end{Bmatrix}$$ are isomorphisms for all irreducible representations $\bS_\nu(V)$ in $\bS_\lambda V\otimes \bS_\mu V$.
 
 Traditionally, motivated by physics for recoupling $3$ angular momenta, the $6j$-symbols were defined for the real Lie group $SU(2)$ whose lie algebra complexification is $\fsl(2)$. Keeping in mind that the convention for physicists is different from that for mathematicians by a factor of $2,$ the irreducibles $V_n$ are parametrized by non-negative half-integers, where $V_n=\Sym^{2n}(\mathbb{C}^2).$ Thus, the parameters for the $6j$-symbols are half-integers, and since the decomposition of the tensor products is multiplicity free by Pieri rule, the symbols are just numbers. Explicit (but complicated) formulas for these numbers are known classically \cite{Varshalovich1989}, and there are calculators and databases available online to calculate them. Using one such calculator \cite{Wigner6jCalculator}, we find that $$ \begin{Bmatrix}
     2 & 1 & 2 \\ 2 & 3 & 2
 \end{Bmatrix}=0.$$ This implies that the flattening map $$\bS_{(4)}V\otimes\bS_{(5,1)}V\to \bS_{(5,1)}V\otimes\bS_{(4)}V$$ is not an isomorphism when $V$ is two dimensional (and for higher dimensions this follows by inheritance). Thus, we see that the flattening maps \eqref{FlatteningMap} and \eqref{FlatteningMap'} are not always isomorphisms. We remark that the formulas for the $6j$-symbols can be used to give exact conditions for when the flattening maps are isomorphisms in the $2$-dimensional case, but since this isn't relevant to our paper given that it is only for two dimensional space, we omit it.  

 While there is a lot of literature on computing $6j$-symbols in low dimensions, see for example \cite{Artamonov_2023} for the $3$ dimensional case, there is very little known in the general case. Theorem \ref{FullrankTheorem} gives a general result about invertibility of certain $6j$-symbols, and we record this as a new theorem. \begin{theorem}
     Suppose $\lambda$ and $\mu$ are two partitions such that $\mu$ is obtained from $\lambda$ by adding $d$ boxes with all boxes added to the same row or column. Then, the $6j$-symbols $$\begin{Bmatrix}
         \lambda & (d) & \mu \\ \lambda & \nu & \mu 
     \end{Bmatrix}\, , \begin{Bmatrix}
         \lambda & (1^d) & \mu \\ \lambda & \nu & \mu 
     \end{Bmatrix}$$ are invertible for all $\bS_\nu V$ in $\bS_\lambda V\otimes \bS_\mu V$ and all dimensions of $V$.
 \end{theorem}

%% file: genericrank.tex
In this section, we give conditions for $\lambda$ and $\mu$ such that the rank of a generic matrix in the image of the map $\phi \colon U\to \mathrm{Hom}(\bS_\lambda V,\bS_\mu V)$ where $U$ is $V$, $\Sym^2 V$ or $\bigwedge^2 V$ is not maximal. We sometimes abuse notation and say that the pair of partitions $(\lambda,\mu)$ is generically non-maximal if this happens. We first discuss the easy case, which is when $U$ is $V$.

\subsection{Generic rank of Pieri Maps over $V$}

We have the following main result.
\begin{proposition}\label{VgeenricRank}
    \begin{enumerate}
        \item All the non-zero matrices in the image of $V$ in $Hom(\bS_{\lambda}V, \bS_{\mu} V)$ have the same rank.
        \item These matrices are of full rank if and only if $\mu$ is obtained from $\lambda$ such that the box is added to the first row or the $n$th row. Further, if the box is added to the first row it is injective and if it is added to the $n$th row it is surjective. 
    \end{enumerate}
\end{proposition}
\begin{proof}
    \begin{enumerate}
        \item This follows from the fact that $\GL(V)$ acts transitively on $V\backslash \{0\}$. 
        
        \item The proof of this fact is given in \cite[Section 2.2]{Wu24}.\qedhere
    \end{enumerate} 
\end{proof}

Theorem \ref{FullrankTheorem} and Proposition \ref{VgeenricRank} together give a complete list of partitions $\lambda$ and $\mu$ (that differ by a box) for which the Young flattening $T_{\lambda \mu}'$ gives a nontrivial lower bound for the border rank of $T_{\lambda\mu}$ - it is exactly when the box is not added to the first or $n$th row.

We were able to solve the problem in this case because $\GL(V)$ acts transitively on $V\backslash \{0\}$, something we do not have for other representations. The idea behind the proof of Proposition \ref{VgeenricRank} was to write $V=\langle v \rangle \oplus H$ and decompose $\bS_\lambda V$ and $\bS_\mu V$ into irreducible representations of $\GL(H)$ and compare the two. It is possible to generalize this argument to study the ranks of matrices that come from an element in the orbit of weight vectors in $\Sym^d(V)$ and $\bigwedge^d V$. But this will not necessarily be the rank of a generic element of $\Sym^d V$ or $\bigwedge^d V$ and so we need to adopt a different strategy to study the generic rank, which is what we discuss in the next subsection.  

\subsection{General Strategy to Study the Generic Rank of $Im(\phi)$}
For notational convenience later, we work with $V^*$ instead of $V$ everywhere. Recall that we use $U$ to denote $\Sym^d V$ or $\bigwedge^d V$. We have the Pieri map $$\bS_\lambda V^*\otimes U^*\to \bS_\mu V^*,$$ which gives us an element in the tensor product $\bS_\lambda V\otimes U\otimes\bS_\mu V^*$. We can view this tensor as a map
 $$\psi\colon U^*\to \mathrm{Hom}(\bS_\mu V,\bS_\lambda V).$$ 
It is easy to see that the matrices in the image under this map are exactly the transposes of the matrices in $\mathrm{Im}(\phi)$. Since the transpose of a matrix has the same rank as the matrix, we can study the generic rank in $\mathrm{Im}(\psi)$ instead. The reason for doing this will become clear soon. 

For the rest of the paper, we restrict our attention to the case when $$\dim(\bS_\mu V)\geq \dim(\bS_\lambda V)$$ and so, generic rank not being maximal is the same as saying that a generic matrix in $\mathrm{Im}(\psi)$ is not surjective. We say that $(\lambda,\mu)$ is not generically surjective when this happens.

\begin{remark}
       We explain why making the assumption that $\dim(\bS_\lambda V)\leq\dim(\bS_\mu V)$ is not that restrictive. From the hook-length formula, we see that 
        \begin{equation}\label{hooklengthequation}
        \frac{\dim(\bS_\mu V)}{\dim(\bS_\lambda V)}=\frac{h_\lambda}{h_\mu}\prod_{(i,j)\in \mu\backslash \lambda}(n-i+j)
        \end{equation}
        where $h_\beta$ denotes the product of the hook lengths for a partition $\beta$. \eqref{hooklengthequation} implies that when $n$ is large enough, the left-hand side is greater than $1$. Hence, for any pair of partitions $(\lambda,\mu)$ with $\lambda\subset \mu$ and $n$ large enough, $\dim (\bS_\mu V)\geq \dim (\bS_\lambda V)$ is satisfied. 
    \end{remark}

Consider the affine space $U^*$ and let $A=\Sym U$ denote the coordinate ring of $U^*$. The Pieri map induces the following map of $A$-modules $$\bS_\mu V\otimes A\to \bS_\lambda V\otimes A(1)\to M\to 0$$ where $M$ is the cokernel. For a fixed element $v\in U^*$, we observe that the restriction of the above map to the fibers over $ v\in U^*$ is exactly the matrix $\psi(v)$. Thus, $\psi(v)$ is surjective if and only if $v$ is not in the support of $M$. In particular, the generic rank in $\mathrm{Im}(\psi)$ is not full if and only if the generic point of $U^*$ is in the support of $M$, which happens exactly when $M$ is supported everywhere. By the definition of rank of an $A$-module $M$, it is clear that the generic point is in the support of $M$ if and only if the rank of $M$ is nonzero. Thus, we have the following proposition. 

\begin{proposition}\label{strategyprop}
    For partitions $\lambda$ and $\mu$ such that $\dim(\bS_\mu V)\geq \dim(\bS_\lambda V)$, the generic rank of matrices in the image of $U\to \mathrm{Hom}(\bS_\lambda V,\bS_\mu V)$ is not maximal if and only if the cokernel of the $A$-linear map  
    \begin{equation}\label{Alinearmap}        
    \bS_\mu V\otimes A\to \bS_\lambda V\otimes A(1)
    \end{equation}
    has positive rank.
\end{proposition}

Thus, our strategy will be to construct resolutions over $A$ that contain the $A$-linear map \eqref{Alinearmap} and use that to determine the rank of the cokernel and then use Proposition \ref{strategyprop} to determine the generic rank.  

\subsection{Generic Rank of Pieri Maps over $\Sym^2 V$}\label{Subsection Sym^2}
For this section, we restrict our attention to partitions $\lambda$ and $\mu$ such that $\mu$ is obtained from $\lambda$ by adding two boxes to the same row. We make this assumption because we know $T_{\lambda,\mu}'$ is full rank only for this case by Theorem \ref{FullrankTheorem}. 

We now explain how we use Kostant's Theorem to construct resolutions over the polynomial ring $A=\Sym (\Sym^2V)$. Let $W=V\oplus V^*$ and consider the symplectic form on $W$ given by: $$\omega((v,\phi),(v',\phi'))=\phi(v')-\phi'(v).$$ The space of endomorphisms of $W$ that fix the symplectic form is the Lie algebra $\fsp(W)$, which is the finite-dimensional simple Lie algebra of type $C_n$.

Using the notation from Theorem \ref{Kostant}, we let $\fg$ be $\fsp(W)$ and take $\fp_Y$ to be $\fp_n$, the parabolic subalgebra corresponding to the $n$th vertex in terms of Bourbaki labelling of the Dynkin diagram. Then, the the Levi subalgebra $\fg_Y$ is given by $\fgl(V)$, which is a reductive Lie algebra with semisimple part isomorphic to the simple Lie algebra of type $A_{n-1}$. The nilpotent radical $ \fu_Y$ as a representation of the Levi is given by $\Sym^2V$. It is abelian and the enveloping algebra $\fU_{\fu_Y}$ is given by the polynomial ring $A$. 

The dominant weights for $\fsp(W)$ are parameterized by partitions with at most $n$ parts. For a partition $\alpha=(\alpha_1,\dots \alpha_{n-1},\alpha_n)$, let $L(\alpha)$ be the corresponding irreducible highest-weight representation of $\fsp(W)$. Since $\Sym^2V\subset \fsp(W)$, $L(\alpha)$ has the structure of an $A$-module. Theorem \ref{Kostant} and Remark \ref{KostantRemark} then give us a way to obtain the minimal free resolution of $L(\alpha)^*$. Before we compute it, we need to talk about the Weyl group of $\fsp(W)$.     

The Weyl group $W$ of $C_n$ is given by the hyperoctahedral group. We now describe the explicit action of the Weyl group on the weights. Let $\beta=(\beta_1,\beta_2\dots,\beta_n)$ be a weight. $W$ is generated by the simple reflections $s_i$ for $1\leq i\leq n-1$ and $\tau$. The action of $s_i$ is given by interchanging $\beta_i$ and $\beta_{i+1}$, and $\tau$ acts by negating $\beta_n$. The Weyl group $W_n$ of the Levi is the symmetric group given by the subgroup generated by the $s_i$. Recall that $W_n'$ denotes the set of elements in $W$ that are the minimal-length representatives of their coset in $W/W_n$. Let us list the elements $w\in W_n'$ for the first few lengths:
\begin{enumerate}
    \item $l(w)=0$. Clearly $\mathrm{id}$ is the only such element.
    \item $l(w)=1$. $\tau$ is the only such element since all the $s_i$ belong to the coset with $\mathrm{id}$.
    \item $l(w)=2$. Again, there is only one such element, which is $s_{n-1}\tau$. Note that $s_i\tau\notin W_n'$ when $i<n-1$ because in this case we have $s_i\tau=\tau s_i$ and thus it belongs to the coset $\tau W_n$. $s_{n-1}$ however does not commute with $\tau$.
    \item $l(w)=3$. In this case we get two elements, $s_{n-2}s_{n-1}\tau$ and $\tau s_{n-1}\tau$. We note that any such element is of the form $as_{n-1}\tau$ where $a$ is some simple reflection. Since $s_i$ commutes with $s_{n-1}$ and $\tau$ for $i<n-2$, a similar argument to the previous one eliminates the possibility of $a$ being $s_i$. 
\end{enumerate}

We remark that we could also deduce the facts above using the fact that $w\in W_n'$ if and only if $w^{-1} \bullet \beta$ is dominant for $\fgl(V)$ when $\beta$ is a partition (recall that a weight is dominant for $\fgl(V)$ if and only if the entries are weakly-decreasing). 

The sum of the fundamental weights $\rho$ in this case is given by $(n,n-1,\dots,1)$. We now have all the tools required to prove the main result of this section. Before we state it, let us first work out an example. 

\begin{example}
    Let $n=4$ and $\alpha=(1,1,1,0)$. The simple reflections are $s_1$, $s_2$, $s_3$ and $\tau$ (using the same notation above). We have $\rho=(4,3,2,1)$. The following table depicts the dotted action on $\alpha$ for those elements in $W_4'$ whose length is less than or equal to $3$. 
       
    \vspace{0.2in}

    \begin{center}
    \begin{tabular}{|c|c|c|c|}
    \hline
        $l(w)$& $w$ & $w^{-1}\bullet\alpha$& $-(w^{-1}\bullet\alpha)^{\mathrm{opp}}$ \\
        \hline
        $0$&$\mathrm{id}$& $(1,1,1,0)$ & $(0,-1,-1,-1)$\\
        $1$& $\tau$ & $(1,1,1,-2)$ & $(2,-1,-1,-1) $\\
        $2$ & $s_3 \tau$ & $(1,1,-1,-4)$ & $(4,1,-1,-1)$ \\
        $3$ & $\tau s_3 \tau$, $s_2 s_3 \tau$ & $(1,1,-3,-4)$,$(1,0,-1,-5)$& $(4,3,-1,-1)$, $(5,1,0,-1)$\\
        \hline
    
    \end{tabular}
    \end{center}
    
        \vspace{0.2in}

From Theorem \ref{Kostant} and Remark \ref{KostantRemark}, we see that the first few terms of the resolution of $L(\alpha)^*$ twisted by $(\det V)$ is given by:
 $$\dots \to (\bS_{(5,4)}V\oplus \bS_{(6,2,1)}V)\otimes A \to \bS_{(5,2)}V\otimes A \to  \bS_{3}V\otimes A \to V\otimes A\to L(\alpha)^*\otimes \det V.$$ 
Now, the dimension of $\bS_{(3)}V$ is $20$. Thus, the kernel of $\bS_{3}V\otimes A \to V\otimes A$ has rank $20-4=16$. It follows that the cokernel of $\bS_{(5,4)}V\oplus \bS_{(6,2,1)}V\otimes A \to  \bS_{(5,2)}V\otimes A$ also has rank $16$. Since the cokernel of the above map is a quotient of the cokernel of the map $\bS_{(5,4)}V\otimes A \to  \bS_{(5,2)}V\otimes A$, it must have rank at least $16$ (it turns out that the rank is still exactly $16$). Thus, this map is not generically surjective. Further, since the dimension of $\bS_{(5,4)}V$ is $280$ and the dimension of $\bS_{(5,2)}V$ is $224$, this proves that the generic rank is not full. 
\end{example}

We are now ready to state the main theorem of this section. 

\begin{theorem}\label{Sym2Pieri}
    Let $\lambda=(\lambda_1,\dots,\lambda_{n})$ and $\mu=(\mu_1,\dots,\mu_n)$ be partitions such that $\mu$ is obtained from $\lambda$ by adding two boxes to the second row. Assume that the following constraint on $\mu$ is satisfied:
   \begin{equation}\label{sym^2Constraint}
       \mu_1+\mu_3+3\leq 2\mu_2.
   \end{equation}
    Then a generic linear map in the image of $\Sym^2 V^*\to \mathrm{Hom}(\bS_{\mu}V, \bS_{\lambda} V)$ is not surjective.
\end{theorem}
\begin{proof}
      Let $\alpha=(\alpha_1,\dots,\alpha_{n-1})$ be a partition with at most $n-1$ parts. Let $\mathbb{F}_\bullet$ denote the minimal free resolution of the module $L(\alpha)^*$. We use Theorem \ref{Kostant} with Remark \ref{KostantRemark} to compute the weights that appear in the first few syzygies of $\mathbb{F}_\bullet$. The weights are displayed in the following table:

    \vspace{0.2in}
    \begin{center}
    \begin{tabular}{|c|c|c|}
    \hline
        $l(w)$& $w$ & $-(w^{-1}\bullet\alpha)^{\mathrm{opp}}$ \\
        \hline
        $0$& $\mathrm{id}$ & $(0,-\alpha_{n-1},\dots,-\alpha_1)$\\
        $1$&$\tau$& $(2,-\alpha_{n-1},\dots,-\alpha_1)$\\  
        $2$ & $s_{n-1} \tau$ & $(\alpha_{n-1}+3,1,-\alpha_{n-2},\dots ,-\alpha_1)$ \\
        $3$ & $\tau s_{n-1} \tau$  & $(\alpha_{n-1}+3, 3,-\alpha_{n-2},\dots ,-\alpha_1)$ \\
        $3$ & $s_{n-2} s_{n-1} \tau$ & $(\alpha_{n-2}+4,1,-\alpha_{n-1}+1,-\alpha_{n-3},\dots,-\alpha_1)$ \\
        \hline
  
    \end{tabular}
      \end{center}
      
    \vspace{0.2in}
    
  We only show the computation for $-(w^{-1}\bullet \alpha)^{\mathrm{opp}}$ when $w=s_{n-1}\tau$ since the computations for the other cases are similar. We have 
  \begin{align*}
      \tau s_{n-1}(\alpha+\rho)-\rho &=\tau s_{n-1}(\alpha_1+n,\dots,\alpha_{n-1}+2,1)-(n,\dots,2,1) \\&=
      \tau (\alpha_1+n,\dots,1,\alpha_{n-1}+2)-(n,\dots,2,1) \\&=
      (\alpha_1+n,\dots,1,-\alpha_{n-1}-2)-(n,\dots,2,1)\\ & =(\alpha_1,\dots,-1,-\alpha_{n-1}-3).
  \end{align*}
   Finally, taking duals, we get $-((s_{n-1}\tau)^{-1}\bullet \alpha)^{\mathrm{opp}}=(\alpha_{n-1}+3,1,-\alpha_{n-2},\dots ,-\alpha_1)$.

   Thus, the first few terms of $\mathbb{F}_\bullet$ twisted by $(\det V)^{\alpha_1}$ (we do this to make all the weights non-negative) is given by: 
    $$\dots \to (\bS_{\mu}V\oplus \bS_{\nu}V)\otimes A \to \bS_{\lambda}V\otimes A \to  \bS_{\alpha'}V\otimes A \to (\det V)^{\alpha_1}\otimes\bS_{\alpha}V^*\otimes A$$ 
    where 

       $$\alpha'=(2+\alpha_1,\alpha_1-\alpha_{n-1},\dots,\alpha_1-\alpha_2),$$ 
        $$\lambda=(\alpha_1+\alpha_{n-1}+3,\alpha_1+1,\alpha_1-\alpha_{n-2},\dots ,\alpha_1-\alpha_2),$$
        $$\mu=(\alpha_1+\alpha_{n-1}+3,\alpha_1+3,\alpha_1-\alpha_{n-2},\dots ,\alpha_1-\alpha_2),$$
       $$\nu=(\alpha_1+\alpha_{n-2}+4,\alpha_1+1,\alpha_1-\alpha_{n-1}+1,\alpha_1-\alpha_{n-3},\dots,\alpha_1-\alpha_2).$$ 
    We observe that the condition $\alpha_{n-2}\geq \alpha_{n-1}$ implies the constraint \eqref{sym^2Constraint}, and conversely if \eqref{sym^2Constraint} is satisfied, then we can find a partition $\alpha$ such that $\mu$ is given in terms of $\alpha$ as above (if $\mu_n\neq 0$, we need to tensor with a suitable power of the determinant, which doesn't change the generic rank). 
    
    Let $K$ denote the cokernel of the map $\mathbb{F}_3\to \mathbb{F}_2$, which coincides with the kernel of the map $\mathbb{F}_1\to \mathbb{F}_0$. $K$ has rank equal to the difference of the ranks of $\mathbb{F}_1$ and $\mathbb{F}_0$. Suppose we know that this is non-zero. Since $K$ is a quotient of the cokernel $M$ of the map \eqref{Alinearmap}, this would imply that $M$ also has positive rank. Thus, we are done by Proposition \ref{strategyprop} once we show that $\dim((\det V)^{\alpha_1}\otimes\bS_{\alpha}V^*)<\dim(\bS_{\alpha'}V)) $.  
    
     We have $(\det V)^{\alpha_1}\otimes\bS_{\alpha}V^*\cong\bS_{\beta}V$ where $\beta=(\alpha_1,\alpha_1-\alpha_{n-1},\dots \alpha_1-\alpha_2)$. We want to show that $$\frac{\dim \bS_{\beta}V}{\dim \bS_{\alpha'}V}<1.$$ We compare the dimensions using the hook-length formula.  We observe that $\alpha'$ is obtained from $\beta$ by adding two boxes to the top row. Thus, the only difference in the numerators of the dimensions according to the hook-length formula is an extra factor of $(n+\alpha_1)(n+\alpha_1+1)$ for $\bS_{\alpha'}V$. Clearly, the hook lengths for the two partitions are the same for all boxes below the first row. Further, we have that $h_{\alpha'}(1,i+2)\leq h_{\beta}(1,i)$ for $1\leq i\leq \alpha_1$ and that $h_{\alpha'}(1,1)\leq \alpha_1+n$ and $h_{\alpha'}(1,2)\leq \alpha_1+n-1$ hold. Combining all these inequalities, we find that $$\frac{\dim \bS_{\beta}V}{\dim \bS_{\alpha'}V}\leq \frac{(\alpha_1+n)(\alpha_1+n-1)}{(\alpha_1+n)(\alpha_1+n+1)}<1. \qedhere$$ 

     \end{proof}

    \begin{remark}
The family of partitions $\lambda$ and $\mu$ from Theorem \ref{Sym2Pieri} has constraints only on the first two rows (which isn't too surprising since the two boxes are added to the second row). A natural question to ask is if these are all the pairs of partitions $(\lambda,\mu)$ (we restrict our attention to when $\mu$ differs from $\lambda$ in the second row only) that are not generically surjective. We have yet to find an example outside of the ones from Theorem \ref{Sym2Pieri} for which $(\lambda,\mu)$ is not generically surjective. 

In particular, if we restrict to the case when we only have two rows, the family of partitions from Theorem \ref{Sym2Pieri} are of the form $\lambda = (\lambda_1,\lambda_2)$ where $$\lambda_1-2\geq\lambda_2\geq  \left\lfloor \frac{\lambda_1}{2}\right\rfloor.$$ We used Macaulay2 to test that for a lot of cases when $\lambda_2<\left\lfloor \frac{\lambda_1}{2}\right\rfloor$ holds, the generic rank is full. It will be interesting to see a proof of this fact if it is indeed true in general.             
    \end{remark}

    \begin{remark}\label{VpieriRemark}
    The pure-free resolutions from \cite{Eisenbud_2011} mentioned in Example \ref{FullrankExample} can be recovered using Kostant's Theorem by considering the case when $\fg$ is the simple Lie algebra $A_n$ and $\fp_Y$ is the parabolic corresponding to the $n$th node. We get a resolution over $\Sym(V)$ in this case, and all the representations in the resolution are irreducible. Further, as mentioned before, for any pair $(\lambda,\mu)$ that differ by a box there exists a pure-free resolution containing the map $$\bS_\mu V\to \bS_\lambda V\otimes \Sym(V)$$ as one of the differentials. If the added box is in the first row, then the map appears as the very first differential and if the box is added to the $n$th row, as the very last differential. Thus, this gives another proof of Proposition \ref{VgeenricRank}.
\end{remark}

    \subsection{Generic Rank of Pieri Maps over $\bigwedge^2 V$}\label{subsection wedge^2}
For this section, we consider partitions $\lambda$ and $\mu$ such that $\mu$ is obtained from $\lambda$ by adding two boxes to the same column. Since the techniques involved are very similar to the techniques in the previous section, we will skip a lot of details and mainly go over the differences from the previous section.  

Let $W=V\oplus V^*$ as in the previous section but this time we work with the symmetric and non-degenerate form $$\omega((v,\phi),(v',\phi'))=\phi(v')+\phi'(v).$$ The space of endomorphisms of $W$ that fix $\omega$ is the Lie algebra $\fso(W)$, which is the finite dimensional simple Lie algebra of type $D_n$.

In this section, we work with $\fso(W)$ and take the parabolic subalgebra $\fp_Y$ to be the one corresponding to the $n$th node (in terms of Bourbaki labelling). The Levi is given by $\fgl(V)$ and the nilpotent radical $\fu_Y$ as a representation of the Levi is $\bigwedge^2V$. The universal enveloping algebra $\fU_{\fu_Y}$ is the polynomial ring $A:=\Sym(\bigwedge^2V)$.

The dominant weights for $\fso(W)$ are parameterized by weakly decreasing tuples $\alpha=(\alpha_1,\dots,\alpha_{n-1},\alpha_n)$ such that either all the $\alpha_i$'s are integers or all of them are half-integers, and $|\alpha_n|\leq \alpha_{n-1}$ holds. The case when all of them are half-integers does not produce any new examples and hence we assume that they are all integers. Let $L(\alpha)$ denote the irreducible finite-dimensional representation corresponding to $\alpha$.

The Weyl group $W$ of $D_n$  is an index $2$ subgroup of the hyperoctahedral group generated by the simple reflections $s_i$ for $1\leq i\leq n-1$ and $\tau$. The $s_i$ act on a weight $\beta=(\beta_1,\dots,\beta_n)$ by interchanging $\beta_i$ and $\beta_{i+1}$. The action of $\tau$ is by negating and interchanging $\beta_{n-1}$ and $\beta_n$. The Weyl group $W_n$ of the Levi is the symmetric group generated by the $s_i$. We list the elements of $W_n'$ with length at most $3$ next.

\vspace{0.2in}

    \begin{center}
    \begin{tabular}{|c|c|}
    \hline
       $l(w)$  & $w$ \\
         \hline
         $0$ & $\mathrm{id}$ \\
         $1$ & $\tau$\\
         $2$ & $s_{n-2}\tau$\\
         $3$ & $s_{n-1}s_{n-2}\tau, s_{n-3}s_{n-2}\tau$\\    
         \hline
    \end{tabular}

\end{center}

\vspace{0.2in}

While we will not give an explanation for the above table since it is similar to the case of $\Sym^2 V$, we will however mention that unlike in the $\Sym^2 V$ case, $\tau s_{n-2}\tau\notin W_n'$ since $\tau s_{n-2}\tau=s_{n-2}\tau s_{n-2}$ (can be verified to be the same signed permutation by explicitly looking at the actions of both elements) and hence shares a coset with $s_{n-2}\tau$. 

Finally, the sum of the fundamental weights is $\rho=(n-1,\dots,0)$. We now state the main Theorem of this section.

\begin{theorem}\label{Wedge2Pieri}
    Let $\lambda=(\lambda_1,\dots,\lambda_n)$ and $\mu=(\mu_1,\dots,\mu_n)$ be partitions such that $\mu$ is obtained from $\lambda$ by adding two boxes to the second and third row in the same column. In particular, we have $\lambda_2=\lambda_3$ and $\mu_2=\mu_3$. Suppose the following constraint holds: 
    \begin{equation}\label{wedge^2constraint}
        \mu_1+\mu_4+2\leq 2\mu_2.
    \end{equation}
     Then a generic linear map in the image of $\bigwedge^2 V^*\to \mathrm{Hom}(\bS_{\mu}V, \bS_{\lambda} V)$ is not surjective.
\end{theorem}

\begin{proof}
     Let $\alpha=(\alpha_1,\dots,\alpha_{n-2})$ be a partition with at most $n-2$ parts. Let $\mathbb{F}_\bullet$ denote the minimal free resolution of the $A$-module $L(\alpha)^*$. Note that $\alpha=(\alpha_1,\dots,\alpha_{n-2},0,0).$ The weights appearing in $\mathbb{F}_i$ for $i\leq 3$ are computed using Theorem \ref{Kostant} and listed in the table below:

    \vspace{0.2in}

    \begin{center}
    \begin{tabular}{|c|c|c|}
    \hline
        $l(w)$& $w$ & $-(w^{-1}\bullet\alpha)^{\mathrm{opp}}$ \\
        \hline
        $0$& $\mathrm{id}$ & $(0,0,-\alpha_{n-2},\dots,-\alpha_1)$\\
        $1$&$\tau$& $(1,1,-\alpha_{n-2}\dots,-\alpha_1)$\\  
        $2$ & $s_{n-2} \tau$ & $(\alpha_{n-2}+2, 1, 1,-\alpha_{n-3},\dots ,-\alpha_1)$ \\
        $3$ & $s_{n-1} s_{n-2}\tau$ & $(\alpha_{n-2}+2, 2, 2,-\alpha_{n-3},\dots ,-\alpha_1)$\\
        $3$ & $s_{n-3} s_{n-2} \tau$ & $(\alpha_{n-3}+3,1, 1,-\alpha_{n-2}+1,-\alpha_{n-4},\dots,-\alpha_1)$ \\
        \hline
    
    \end{tabular}
    \end{center}

    \vspace{0.2in}

    After twisting by a power of the determinant, the first few terms of $\mathbb{F}_\bullet$ look like: $$\dots \to (\bS_{\mu}V\oplus \bS_{\nu}V)\otimes A \to \bS_{\lambda}V\otimes A \to  \bS_{\alpha'}V\otimes A \to (\det V)^{\alpha_1}\otimes\bS_{\alpha}V^*\otimes A$$  where 

       $$\alpha'=(\alpha_1+1,\alpha_1+1, \alpha_1-\alpha_{n-2},\dots,\alpha_1-\alpha_2),$$ 
       $$\lambda=(\alpha_1+\alpha_{n-2}+2,\alpha_1+1,\alpha_1+1,\alpha_1-\alpha_{n-3},\dots ,\alpha_1-\alpha_2),$$
       $$\mu=(\alpha_1+\alpha_{n-2}+2,\alpha_1+2,\alpha_1+2,\alpha_1-\alpha_{n-3},\dots ,\alpha_1-\alpha_2),$$
       $$\nu=(\alpha_1+\alpha_{n-3}+3,\alpha_1+1,\alpha_1+1, \alpha_1-\alpha_{n-2}+1,\alpha_1-\alpha_{n-4},\dots,\alpha_1-\alpha_2).$$ 
We observe that $\alpha_{n-3}\geq \alpha_{n-2}$ is equivalent to the inequality \eqref{wedge^2constraint} on $\mu$.

The rest of the argument is almost the same as the argument for the $\Sym^2 V$ case. The only thing to check is if the rank of $\mathbb{F}_1$ is strictly greater than the rank of $\mathbb{F}_0$ in this case. This can be shown using the hook-length formula and the calculation is very similar to the $\Sym^2 V$ case, and hence we skip it.    
\end{proof}

\begin{remark}
    We note that we restricted our attention to the case when the two boxes are added to the same row or column since that is when we know that the Young Flattening $T_{\lambda,\mu}'$ is full rank. However, the techniques in the proofs of Theorems \ref{Sym2Pieri} and \ref{Wedge2Pieri} can be used to produce other families of partitions $\lambda$ and $\mu$ where the added boxes are not necessarily in the same row or column by looking at the minimal free resolutions of $L(\alpha)$ with different constraints on $\alpha$. We list the families below along with the constraints on $\alpha$ that produce these families.

\begin{enumerate}
    \item $U=\Sym^2 V$, $\alpha$ dominant weight of $C_n$
    
    \begin{enumerate}
        \item $\alpha_{n-1}=\alpha_n$
        $$\lambda=(\alpha_1+\alpha_{n-1}+2,\alpha_1-\alpha_{n-1},\alpha_1-\alpha_{n-2},\dots,\alpha_1-\alpha_2),$$
        $$\mu=(\alpha_1+\alpha_{n-1}+3,\alpha_1-\alpha_{n-1}+1,\alpha_1-\alpha_{n-2},\dots,\alpha_1-\alpha_2).$$
        \item $\alpha_{n-2}=\alpha_{n-1}$
        $$\lambda=(\alpha_1+\alpha_{n-2}+3,\alpha_1-\alpha_n+1,\alpha_1-\alpha_{n-2},\dots,\alpha_1-\alpha_2),$$
         $$\mu=(\alpha_1+\alpha_{n-2}+4,\alpha_1-\alpha_n+1,\alpha_1-\alpha_{n-2}+1,\dots,\alpha_1-\alpha_2).$$
    \end{enumerate}

    \item $U=\bigwedge^2 V$, $\alpha$ integer-dominant weight of $D_n$

    \begin{enumerate}
        \item $\alpha_{n-1}=\alpha_n$
        $$\lambda=(\alpha_1+\alpha_{n-2}+2,\alpha_1+\alpha_{n-1}+1,\alpha_1-\alpha_{n-1}+1,\alpha_1-\alpha_{n-3},\dots,\alpha_1-\alpha_2),$$
        $$\mu=(\alpha_1+\alpha_{n-2}+2,\alpha_1+\alpha_{n-1}+2,\alpha_1-\alpha_{n-1}+2,\alpha_1-\alpha_{n-3},\dots,\alpha_1-\alpha_2).$$
        \item $\alpha_{n-2}=\alpha_{n-1}$
         $$\lambda=(\alpha_1+\alpha_{n-2}+1,\alpha_1+\alpha_n+1,\alpha_1-\alpha_{n-2},\dots,\alpha_1-\alpha_2),$$
         $$\mu=(\alpha_1+\alpha_{n-2}+2,\alpha_1+\alpha_n+1,\alpha_1-\alpha_{n-2}+1,\dots,\alpha_1-\alpha_2).$$
         \item $\alpha_{n-3}=\alpha_{n-2}$
         $$\lambda=(\alpha_1+\alpha_{n-3}+2,\alpha_1+\alpha_n+1,\alpha_1-\alpha_{n-1}+1,\alpha_1-\alpha_{n-3},\dots,\alpha_1-\alpha_2),$$
         $$\mu=(\alpha_1+\alpha_{n-3}+3,\alpha_1+\alpha_n+1,\alpha_1-\alpha_{n-1}+1,\alpha_1-\alpha_{n-3}+1,\dots,\alpha_1-\alpha_2).$$
         
    \end{enumerate}
\end{enumerate}

    We remark that the $\lambda$ and $\mu$ in the families we obtain all differ only in the first few rows. The reason for this is because we only look at the first few syzygies of the minimal free resolution from Kostant's theorem since this is where the free modules are irreducible $\fgl(V)$ representations. In the next section we will see how to extend the family of examples to the cases where the boxes are added not just to the top few rows.       
\end{remark}

\subsection{Extending the Families of Partitions with non-maximal Generic Rank}
In Sections \ref{Subsection Sym^2} and \ref{subsection wedge^2} we constructed examples where the two boxes were added to the second row in the $\Sym^2 V$ case and to the second and third rows (same column) in the $\bigwedge^2 V$ case. In the $\Sym^2 V$ case, if the two boxes are added to the first row, and in the $\bigwedge^2 V$ case if the boxes are added to the first and second rows, the corresponding tensor $T_{\lambda \mu}$ is of minimal border rank (\cite[Proposition 1]{Wu24}). In this section, we will see how to construct families of pairs of partitions $(\lambda,\mu)$ that are generically non-maximal when the boxes are not necessarily added to the first or second row. 

Before we state the main theorem of this section, we first introduce some notation. For a partition $\lambda$ and an integer $k\geq \lambda_1$, let $\lambda^k$ denote the new partition $$\lambda^k=(k,\lambda_1,\lambda_2,\dots).$$ 

\begin{theorem}\label{PieriBWBTheorem}
     Let $\lambda$ and $\mu$ be partitions such that $\mu$ is obtained from $\lambda$ by adding $d$ boxes with no two boxes added to the same column (or same row in the $\bigwedge^d$ case). For a vector space $W$ of dimension $n-1$, let $B$ denote $\Sym U'$ where $U'$ is $\Sym^d W$ or $\bigwedge^d W$ depending on context. Suppose we know that the cokernel $M$ of $$\bS_\mu W\otimes B\to \bS_\lambda W\otimes B(1)$$ has positive rank. Then the same holds for partitions $\lambda^k$ and $\mu^k$ for all but finitely many integers $k$ over a vector space $V$ of dimension $n$. 
 \end{theorem}

 For the $\Sym^2 V$ and $\bigwedge^2 V$ cases, Theorem \ref{PieriBWBTheorem} along with Proposition \ref{strategyprop} implies that we have infinitely many family of partitions with non-maximal generic rank when the boxes are added to any row below the first row (in the same row/column) by starting with the family of examples from Theorems \ref{Sym2Pieri} and \ref{Wedge2Pieri} and using the above construction iteratively.

 Before we give the proof of Theorem \ref{PieriBWBTheorem}, we sketch the idea of the proof by working out an example. 

 \begin{example}
     Consider the case when $\lambda=(3,1)$ and $\mu=(3,3)$ and $W$ is of dimension $4$. Theorem \ref{Sym2Pieri} says that the cokernel $M$ in this case has positive rank. In fact, its resolution is given by: 

     $$ 0\to \bS_{(5^4)}W\to \bS_{(5^3,3)} W\to \bS_{(5^2,4,2)}W\to \bS_{(5^2,2^2)}W\to \bS_{(3,3)}W\to \bS_{(3,1)}W$$
where $\bS_\alpha W$ denotes the free $B$-module $\bS_\alpha W\otimes B$. Working over $X=\Gr(4,5)$, we sheafify the above complex by replacing $W$ with the rank $4$ tautological bundle $\cR$ over $X$ and base change along the flat extension $\Sym(\Sym^2 \cR)\to \Sym (\Sym^2 V)$. We use $A$ to denote $\Sym(\Sym^2 V)$. We then tensor the resulting acyclic complex with $\cO_X(k)$ for $k>5$ and take global sections. By Borel-Weil-Bott, each of the terms in the complex has non-vanishing cohomologies in degree $0$ only, and we have $$H^0(X,\bS_\alpha\cR(k)\otimes A)=\bS_{\alpha^k}V\otimes A.$$ Hence, we obtain the following acyclic complex $H_\bullet(k)$: $$ 0\to \bS_{(k,5^4)}V\to \bS_{(k,5^3,3)} V\to \bS_{(k,5^2,4,2)}V\to \bS_{(k,5^2,2^2)}V\to \bS_{(k,3,3)}V\to \bS_{(k,3,1)}V$$ where we use $\bS_\alpha V$ to denote $\bS_\alpha V\otimes A$. 

Using the hook-length formula to compute the dimensions (as polynomials in $k$) of the Schur modules above, we obtain the Euler characteristic of $H_\bullet(k)$ for $k>5$: $$\frac{(k+2)(k+3)(k+4)(9k-11)}{24}.$$ This polynomial has no positive integer roots and hence the Euler characteristic is non-zero for all $k>5$. This implies that the cokernel $$\bS_{(k,3,3)}V\otimes A\to \bS_{(k,3,1)} V\otimes A(1)\to M_k\to 0$$ has positive rank for $k>5$ and hence by Proposition \ref{strategyprop} we have constructed a family of partitions that are not generically surjective, where the boxes are added to the third row.
 \end{example}
 
\subsection*{Proof of Theorem \ref{PieriBWBTheorem}.} 
    Let $F_\bullet$ denote the minimal free resolution of $M$ as a $B$-module:
    $$\dots\to F_i\to\dots \to F_2\to \bS_\mu W\otimes B(-1)\to \bS_\lambda W\otimes B\to M\to 0.$$
    Since this is an exact sequence, the Euler characteristic of $F_\bullet$ is equal to zero. If $G_\bullet$ is the complex obtained from $F_\bullet$ by excluding $M$, then the Euler characteristic of $G_\bullet$ is equal to the rank of $M$ and hence is non-zero. 

    Let $X=\Gr(n-1,n)$ be the Grassmannian consisting of hyperplanes inside the $n$ dimensional complex vector space $V$. We have the following tautological exact sequence of vector bundles over $X$ 
$$0\to \cR\to V\otimes \cO_X\to \cO_X(1)\to 0.$$
The rank of $\cR$ equals $n-1$. Hence, we can sheafify the complex $G_\bullet$ by replacing $W$ with $\cR$ to obtain an acyclic complex of vector bundles over $X$. We base change from $\Sym (\cU')$ to $\Sym(U)=A$ where $\cU'$ denotes $\Sym^d \cR$ or $\bigwedge^d \cR$ based on context. We note that the inclusion $\Sym(\cU')\to A$ is flat (locally it is an inclusion of polynomial rings) and hence this base-change preserves exactness. We tensor this complex throughout with $\cO_X(k)$. Let $l$ be the largest size of the first row of any partition that appears as a weight in the complex $G_\bullet$. Then for $k>l$ and any partition $\nu$ that appears in some $G_i$, we have by Borel-Weil-Bott that $$H^i(X,\bS_\nu \cR(k)\otimes A)=\begin{cases}
    \bS_{\nu^k}V \otimes A& \text{when $i=0$}\\ 0 & \text{otherwise}
\end{cases}$$
since $(k,\nu_1,\dots,\nu_{{n-1}})$ is already weakly decreasing and $A$ is a trivial vector bundle over $X$. Thus when $k>l$, from the $E_1$ page of the hypercohomology spectral sequence we obtain the complex of $A$-modules with each weight $\nu$ in $G_\bullet$ replaced by $\nu^k$ and $W$ replaced by $V$ everywhere. We denote this complex by $H_\bullet(k)$. The first differential of $H_\bullet(k)$ is the map $$\bS_{\mu^k}V\otimes A(-1)\to \bS_{\lambda^k} V\otimes A$$ whose cokernel we want to show has non-zero rank. By running the spectral sequence horizontally first, we see that $H_\bullet(k)$ is acyclic. Thus, it suffices to show that the Euler characteristic of $H_\bullet(k)$ is non-zero for all but finitely many integers $k>l$. 

For a partition $\nu$, let us compare the dimensions of $\bS_\nu W$ and $\bS_{\nu^k} V$. We use hook length formula. 

We have $$\dim(\bS_{\nu^k}V)=(\prod_{i>1,j} \frac{n-i+j}{h_{\nu^k}(i,j)})(\prod_{j=1}^{k}\frac{n-1+j}{h_{\nu}(1,j)}).$$
For the first term, we can make the substitution $i=i'+1$ from which we get 
$$\prod_{i>1,j} \frac{n-i+j}{h_{\nu^k}(i,j)}=\prod_{i',j} \frac{n-1-i'+j}{h_{\nu}(i',j)}=\dim(\bS_\nu W).$$
The numerator of the second term can be evaluated as follows:
$$\prod_{j=1}^{k}(n-1+j)=\frac{(n-1+k)!}{(n-1)!}.$$
We observe that $$k+n-1\geq h_{\nu}(1,1)>h_\nu(1,2)>\dots h_\nu(1,k)\geq1.$$ 
Thus, $$\prod_{j=1}^{k}\frac{n-1+j}{h_{\nu}(1,j)}=\frac{(n-1+k)!}{(n-1)!(\prod_{j=1}^kh_\nu(1,j))}$$ is a degree $n-1$ polynomial in $k$ with leading coefficient $\frac{1}{(n-1)!}$. Combining the two terms, we see that the function $k\to \dim(\bS_{\nu^k}V)$ is a polynomial function of degree $n-1$ with leading coefficient $\frac{\dim(\bS_\nu W)}{(n-1)!}$ for $k>\nu_1$. 

From the above calculation, it follows that the function that sends $k$ to the Euler Characteristic of $H_\bullet(k)$ is a polynomial of degree $n-1$ with leading coefficient equal to the Euler Characteristic of $G_\bullet$ divided by a factor of $(n-1)!$ when $k>l$ holds. Thus, this is a non-zero polynomial and since the roots are finite, it follows that the Euler Characteristic of $H_\bullet(k)$ is non-zero for all but finitely many $k$. \qed

%% file: generalisations.tex
\subsection{Lower Bounds for Border Rank of Tensors Induced by Pieri Maps}\label{Explicit Lower Bound section}
We prove the following theorem using results from previous sections. 

\begin{theorem}\label{ExplicitBoundTheorem}
For partitions $\lambda=(\lambda_1,\dots,\lambda_{n})$ and $\mu=(\mu_1,\dots,\mu_n)$, let $k=\dim \bS_\lambda V$ and $\ell=\dim \bS_\mu V$.  
\begin{enumerate}
    \item Let $U=V$. Suppose $\mu$ is obtained from $\lambda$ by adding a box to the $t$th row, where $t\neq 1,n.$ Then, $T_{\lambda \mu}$ is not of minimal border rank. Further, let 
    $$\lambda'=(\lambda_1,\dots, \lambda_{t-2},\lambda_t,\lambda_t,\dots,\lambda_n),$$  $$\lambda''=(\lambda_1,\dots,\lambda_{t-3},\lambda_{t-1}-1,\lambda_t,\lambda_t,\dots,\lambda_n),$$
  $$\mu'=(\mu_1,\dots,\mu_t,\mu_t,\mu_{t+2},\dots,\mu_n),$$
  $$\mu''=(\mu_1,\dots,\mu_t,\mu_t,\mu_{t+1}+1,\mu_{t+3},\dots,\mu_n).$$
  Let $c_1=\dim \bS_{\lambda'}V - \dim \bS_{\lambda''} V$ and $c_2=\dim \bS_{\mu'} V-\dim \bS_{\mu''} V.$ Then, we have 
  $$\underline{\mathrm{Rk}}(T_{\lambda \mu})\geq \max\{\ell(\frac{k}{k-c_1}),  k(\frac{\ell}{\ell-c_2})\}.$$
    \item Let $U=\Sym^2 V$. Assume that $\ell\geq k$ holds. Suppose $\mu$ is obtained from $\lambda$ by adding two boxes to the second row, and $\mu$ satisfies 
       $$\mu_1+\mu_3+3\leq 2\mu_2.$$
       Then, $T_{\lambda \mu}$ is not of minimal border rank. Further, let $$\alpha'=(\mu_2-1,2\mu_2-\mu_1-3,\mu_3,\dots,\mu_n),$$
$$ \alpha=(\mu_2-3,2\mu_2-\mu_1-3,\mu_3,\dots,\mu_n).$$  Let $c=\dim \bS_{\alpha'} V-\dim \bS_{\alpha} V$. Then we have $$\underline{\mathrm{Rk}}(T_{\lambda \mu})\geq  \ell(\frac{k}{k-c}).$$  
\item Let $U=\bigwedge^2 V$. Assume that $\ell\geq k$ holds. Suppose $\mu$ is obtained from $\lambda$ by adding two boxes to the second and third rows but in the same column, and $\mu$ satisfies 
       $$\mu_1+\mu_4+2\leq 2\mu_2.$$
       Then, $T_{\lambda \mu}$ is not of minimal border rank. Further, let $$\alpha'=(\mu_2-1,\mu_2-1, 2\mu_2-\mu_1-2,\mu_4,\dots,\mu_n),$$
$$ \alpha=(\mu_2-2,\mu_2-2, 2\mu_2-\mu_1-2,\mu_4,\dots,\mu_n).$$  Let $c=\dim \bS_{\alpha'} V-\dim \bS_{\alpha} V$. Then we have $$\underline{\mathrm{Rk}}(T_{\lambda \mu})\geq \ell(\frac{k}{k-c}).$$  
\end{enumerate}
\end{theorem}
\begin{proof}
 Theorem \ref{FullrankTheorem} along with \eqref{k-lLowerBound} gives us the following lower bound: 
 $$\underline{\mathrm{Rk}}(T_{\lambda \mu})\geq \left \lceil\frac{k\ell}{r}\right \rceil,$$ where $r$ is the generic rank in the image of $U$ in $\mathrm{Hom}(\bS_{\lambda}V,\bS_{\mu}V)$. It follows from Theorems \ref{Sym2Pieri} and \ref{Wedge2Pieri} and Remark \ref{VpieriRemark} that for partitions $\lambda$ and $\mu$ as in the Theorem, $r<k,\ell$ holds. This proves that the the tensor $T_{\lambda \mu}$ is not of minimal rank.

 To obtain the explicit bounds, if follows from the proofs of Theorems \ref{Sym2Pieri} and \ref{Wedge2Pieri} that the generic rank $r$ in the $\Sym^2 V$ and $\bigwedge^2 V$ cases satisfies $r\leq k-c,$ since the respective cokernels of \eqref{Alinearmap} have rank greater than or equal to the rank of the cokernel of the third differential in the resolution coming from Kostant's theorem, whose rank is $c$. For the case when $U=V$, we observe that the partitions $\lambda'$, $\lambda''$, $\mu'$ and $\mu''$ are the ones that appear in the pure free resolution constructed in \cite{Eisenbud_2011} that contain the Pieri map between $\lambda$ and $\mu$ as a differential. Thus, the lower bound for this case follows from Remark \ref{VpieriRemark}.  
\end{proof}

\begin{remark}\label{VcaseExplicitBound}
    The lower bound stated in Theorem \ref{ExplicitBoundTheorem} for the case $U=V$ is not the sharpest bound one obtains using the resolution in \cite{Eisenbud_2011}, it is just the bound that is simpler to compute. However, the sharpest lower bounds are obtained by replacing $c_1$ and $c_2$ with $r_1$ and $r_2$, where $r_1$ and $r_2$ are the ranks of the cokernel and kernel of the differential $d$ between the $\bS_\mu V$ and $\bS_\lambda V$ terms in the resolution. This in turn corresponds to the ranks of the differentials adjacent to $d$, which are given by alternating sums of Betti numbers of the resolution. These Betti numbers are equal to the dimensions of certain Schur functors that can be computed using the hook length formula. 

    For example, let $\lambda=(1^r)$ and $\mu=(1^{r+1})$ where $r<n-1$ holds. Then, the pure free resolution from \cite{Eisenbud_2011} containing the required differential is just the Koszul complex $$\dots\to \bigwedge^{r+2}V\otimes A(-1)\to \bigwedge^{r+1} V\otimes A\xrightarrow{*} \bigwedge^{r}V\otimes A(1)\to \bigwedge^{r-1} V\otimes A(2)\to \dots$$ whose Betti numbers are given by the formula $$\beta_i=\binom{n}{i}.$$ Thus, the ranks $r_1$ and $r_2$ of the cokernel and kernel of the starred differential are given by $$r_1=\sum_{i=0}^{r-1}(-1)^{r-1-i} \binom{n}{i}=\binom{n-1}{r-1},$$ $$r_2=\sum_{i=r+2}^{n}(-1)^{r+2-i} \binom{n}{i}=\binom{n-1}{r+1}.$$ We used Pascal's identity for the final equalities in the expressions above. Thus, we have $$\ell-r_2=k-r_1=\binom{n}{r}-\binom{n-1}{r-1}=\binom{n-1}{r}.$$ This gives the following lower bound for border rank:

    \begin{equation}\label{wedge^rExplicitBound}\underline{\mathrm{Rk}}(T_{\lambda\mu})\geq \frac{\binom{n}{r+1}\binom{n}{r}}{\binom{n-1}{r}}=\binom{n}{r+1}\frac{n}{n-r}.\end{equation} 
\end{remark}

We now use the explicit lower bounds obtained to see how well we can do. Consider the example from remark \ref{VcaseExplicitBound} with $\lambda=(1^r)$ and $\mu=(1^{r+1}).$ Let $k$ and $\ell$ denote the dimensions of $\bS_\lambda V$ and $\bS_\mu V$ respectively. Let $r=\left \lceil \frac{n-1}{2} \right \rceil-1$, then we have $\ell\geq k$, and \eqref{wedge^2constraint} gives us the bound $$\underline{\mathrm{Rk}}(T_{\lambda\mu})\geq \ell(2-\frac{t}{r+t})$$ where $t=1$ when $n$ is odd and $t=2$ when $n$ is even. As $n$ goes to infinity, $r$ goes to infinity and the lower bound tends to $2\ell.$ Thus, the lower bound for border rank for the family of invariant tensors in $V\otimes \bigwedge^r V\otimes \bigwedge^{r+1} V$ when $r=\left \lceil \frac{n-1}{2} \right \rceil-1$ induced by the Pieri map $$V\otimes \bigwedge^rV\to \bigwedge^{r+1}V$$ tends to $2\ell$ where $\ell$ is the dimension of the largest ambient vector space. We remark that the Young flattening used for these tensors is the Koszul flattening, which was used by Landsberg and Ottaviani in \cite{v011a011} to give lower bounds for the matrix multiplication tensor.

Next, we show that the Koszul flattening gives the best lower bounds out of all the Young flattenings we study in the case when $U=V.$

\begin{proposition}\label{2lCannotBreachProp}
    Let $\lambda$ and $\mu$ be two partitions with $\mu$ obtained from $\lambda$ by adding a box. Then, the lower bound on border rank for $T_{\lambda \mu}\in V\otimes \bS_\lambda V\otimes \bS_\mu V$ obtained using the Young flattening \eqref{FlatteningMap} is less than $2\ell,$ where $\ell$ is the dimension of the largest ambient vector space. 
\end{proposition}
\begin{proof}
    Assume $\dim \bS_\mu V\geq \dim \bS_\lambda V$ holds. Since the proof for the opposite case is very similar, we only prove it for this case. Let $\ell$ and $k$ denote the dimensions of $\bS_\mu V$ and $\bS_\lambda V$ respectively and suppose that $\mu$ is obtained from $\lambda$ by adding a box to the $t$th row. 
    
    Let $r$ denote the rank of a non-zero matrix in the image $V\to \mathrm{Hom}(\bS_\lambda V, \bS_\mu V)$ induced by the Pieri map. Then, by \ref{k-lLowerBound}, the Young flattening gives us the following lower bound: $$\underline{\mathrm{Rk}}(T_{\lambda \mu})\geq   \ell(\frac{k}{r}).$$ Thus, our goal is to show that $\frac{k}{r}<2$ holds. 
    
    Fix $v\in V$ and let $\phi_v\in \mathrm{Hom}(\bS_\lambda V,\bS_\mu V)$ be the corresponding matrix. Following \cite[Section 2.2]{Wu24}, let $H$ denote an $n-1$ dimensional subspace of $V$ such that $V=\langle v \rangle \oplus H.$ Then, $\bS_\lambda V$ and $\bS_\mu V$ can be decomposed into irreducible $\GL(H)$-representations using branching rule: $$\bS_\lambda V= \bigoplus_{ \nu} \bS_\nu H,$$ $$\bS_\mu V= \bigoplus_{ \nu} \bS_\nu H, $$ where the direct sum is over all $\nu$ obtained from $\lambda$ (resp. $\mu$) by removing at most one box per column. With respect to this decomposition, in \cite{Wu24}, the author shows that the kernel of $\phi_v$ contains the direct sum of all $\bS_\nu H\subset \bS_\lambda V$ for which $\nu_{t-1}=\lambda_t,$ and the rest map injectively to $\bS_\mu V.$ Let $A$ denote the set of all partitions $\nu$ that appear in the decomposition of $\bS_\lambda V$ into $\GL(H)$-irreducibles. Define subsets $B$ and $C$ as follows: $$B=\{\nu \in A \mid \nu_{t-1}=\lambda_t\}, \qquad C=\{\nu\in A\mid \nu_{t-1}>\lambda_t\}.$$ Thus, $B$ is the subset of partitions in the kernel and $A=B\cup C.$ Using this notation, to show $r>\frac{k}{2}$ holds, it suffices to show $$r=\dim \bigoplus_C \bS_\nu H > \dim \bigoplus_B \bS_\nu H=k-r.$$ Consider the map $\psi\colon B\to C$ that corresponds to the operation of adding a box to the $t-1$th row. We will now show that $\dim \bS_\nu H < \dim \bS_{\psi(\nu)}H$ holds. 

    By the hook length formula, we have $$\frac{\dim \bS_{\psi(\nu)}H}{\dim \bS_\nu H}=(n-1+\nu_{t-1}+1-(t-1))\frac{h_\nu}{h_{\psi(\nu)}}$$ where $h_\pi$ denotes the product of the hook lengths for a partition $\pi.$ Now, \begin{align*} \frac{h_\nu}{h_{\psi(\nu)}}&=\frac{\prod_{i=1}^{t-2} h_\nu(i,\nu_{t-1}+1)}{\prod_{i=1}^{t-2} h_{\psi(\nu)}(i,\nu_{t-1}+1)} .\frac{\prod_{i=1}^{\nu_{t-1}} h_\nu(t-1,i)}{\prod_{i=1}^{\nu_{t-1}} h_{\psi(\nu)}(t-1,i)}\\ &= \frac{\prod_{i=1}^{t-2} h_\nu(i,\nu_{t-1}+1)}{\prod_{i=1}^{t-2} (h_{\nu}(i,\nu_{t-1}+1)+1)} .\frac{\prod_{i=1}^{\nu_{t-1}} h_\nu(t-1,i)}{\prod_{i=1}^{\nu_{t-1}} (h_{\nu}(t-1,i)+1)}  \end{align*}  since the other terms cancel out. Now, we know that the function $x\to \frac{x}{x+1}$ is an increasing function. Since $\nu$ is obtained from $\lambda$ by removing at most one box from each column, and keeping in mind that $\nu_{t-1}=\lambda_t,$ we have $ h_\nu(t-1,i)\geq h_\lambda(t,i)$ and $\nu_{j-1}\geq \lambda_j$ to hold, which in turn implies $h_\nu(i,\nu_{t-1}+1)\geq h_\lambda(i+1,\lambda_t+1).$ Combining all this together gives $$ \frac{\prod_{i=1}^{t-2} h_\nu(i,\nu_{t-1}+1)}{\prod_{i=1}^{t-2} (h_{\nu}(i,\nu_{t-1}+1)+1)} .\frac{\prod_{i=1}^{\nu_{t-1}} h_\nu(t-1,i)}{\prod_{i=1}^{\nu_{t-1}} (h_{\nu}(t-1,i)+1)}>  \frac{\prod_{i=1}^{t-1} h_\lambda(i,\lambda_t+1)}{\prod_{i=1}^{t-1} (h_{\lambda}(i,\lambda_t+1)+1)} .\frac{\prod_{i=1}^{\lambda_{t}} h_\lambda(t,i)}{\prod_{i=1}^{\lambda_{t}} (h_{\lambda}(t,i)+1)} .$$ We note that the inequality is strict because the right hand side has the extra  $$\frac{ h_\lambda(1,\lambda_t+1)}{h_{\lambda}(1,\lambda_t+1)+1)}<1$$ term. We observe that the right hand side of the above inequality is just $\frac{h_\lambda}{h_\mu}$ and hence we have established that $$ \frac{h_\nu}{h_{\psi(\nu)}}>\frac{h_\lambda}{h_\mu}$$ holds. Finally, we have the following sequence of inequalities: $$\frac{\dim \bS_{\psi(\nu)}H}{\dim \bS_\nu H}=(n+\nu_{t-1}+1-t)\frac{h_\nu}{h_{\psi(\nu)}}>(n+\lambda_t+1-t)\frac{h_\lambda}{h_\mu}=\frac{\ell}{k}\geq 1. $$ Thus, we have shown that $\dim \bS_{\psi(\nu)}H>\dim \bS_\nu H$ holds, and hence $$\dim \bigoplus_C \bS_\nu H > \dim \bigoplus_B \bS_\nu H$$ holds as required.
\end{proof}

We next present an explicit calculation of the lower bounds obtained in the $\Sym^2 V$ case, establishing the lower bound mentioned in Table \ref{tab:my_label}.

\begin{corollary}
    Let $\lambda=(d,d-2)$, $\mu=(d,d)$ and $n\sim (\sqrt 3-1)d$ ($\frac{n}{d}$ slightly greater than $\sqrt 3-1$). Then the tensor $T_{\lambda\mu}\in \Sym^2V \otimes \bS_\lambda V\otimes \bS_\mu V^*$ has border rank lower bound given by $1.42\ell$, where $\ell=\dim \bS_\mu V.$
\end{corollary}
\begin{proof}
    We use the hook-length formula and cancel out common terms to obtain the ratio of $\dim \bS_\lambda V$ and $\dim \bS_\mu V:$ $$\frac{\dim \bS_\mu V}{\dim \bS_\lambda V}=\frac{(n+d-2)(n+d-3)}{3d(d-1)}.$$ The relation between $n$ and $d$ ensures that $n+d>\sqrt{3}d$ holds and hence, $\dim \bS_\mu V\geq \dim \bS_\lambda V$ holds.

    The partitions $\alpha'$ and $\alpha$ from Theorem \ref{ExplicitBoundTheorem} correspond to $(d-1,d-3)$ and $(d-3,d-3)$ in our case. Using the hook-length formula and simplifying by cancelling common terms, and plugging in the formula for $n$ in terms of $d$, we obtain the following ratios:
    $$\frac{\dim \bS_{\alpha'} V}{\dim \bS_\lambda V}=\frac{(d-2)(d+1)}{(n+d-4)(n+d-1)}\sim \frac{1}{3},$$ $$\frac{\dim \bS_{\alpha} V}{\dim \bS_\lambda V}=\frac{(d-2)(d-1)d(d+1)}{3(n+d-4)(n+d-3)(n+d-2)(n+d-1)}\sim \frac{1}{27}.$$ Plugging these ratios into the lower bound in Theorem \ref{ExplicitBoundTheorem}, we obtain the lower bound $$\underline{\mathrm{Rk}}(T_{\lambda\mu})\geq \ell(\frac{1}{1-\frac{1}{3}+\frac{1}{27}})\sim 1.42\ell. $$ This finishes the proof. 
\end{proof}

We are unaware whether the $2\ell$ bound can be breached when we generalize to $U=\Sym^2 V$ or $U=\bigwedge^2 V.$ The examples coming from Kostant's thoerem are not strong enough to breach the bound since the boxes are constrained to being added to the second row (or second and third rows in the $\bigwedge^2V$ case). However, it is possible that the bound is breached when the boxes are added to lower rows, and towards this direction, we state the following result that is implied by Theorem \ref{PieriBWBTheorem}. 

\begin{theorem}
    Assume that $W$ is an $n-1$ dimensional complex vector space and that $\ell=\dim \bS_\mu W\geq \dim \bS_\lambda W$ holds. Let $k\geq \mu_1$, $\lambda^k=(k,\lambda_1,\dots\lambda_n)$ and $\mu^k=(k,\mu_1,\dots,\mu_n)$. Suppose that the border rank lower bound for $T_{\lambda \mu}\in \Sym^d W\otimes \bS_\lambda W\otimes \bS_\mu W$ (resp. $T_{\lambda \mu}\in \bigwedge^d W\otimes \bS_\lambda W\otimes \bS_\mu W$) obtained using the Young flattening \ref{FlatteningMap} (resp. \ref{FlatteningMap'}) is $\ell C$. Then, the border rank for the family of tensors $T_{\lambda^k \mu^k}\in \Sym^d V\otimes \bS_{\lambda^k}V\otimes \bS_{\mu^k}V$ (resp. $T_{\lambda^k \mu^k}\in \bigwedge^d V\otimes \bS_{\lambda^k}V\otimes \bS_{\mu^k}V$) satisfies the following inequality: 
    $$\liminf_{k\to \infty } \underline{\mathrm{Rk}}(T_{\lambda^k \mu^k})\geq \ell_kC$$ where $\ell_k=\dim \bS_{\mu^k}V.$
\end{theorem}

\begin{proof}
 We know from the explicit bounds in Theorem \ref{ExplicitBoundTheorem} that the factor $$C=\frac{\dim \bS_\lambda W}{\dim \bS_\lambda W-\mathrm{rk} M}$$ where $M$ is the cokernel of \eqref{Alinearmap}. Thus, it suffices to show that $$\lim_{k\to \infty} \frac{\mathrm{rk} M_k}{\dim \bS_{\lambda^k }V}=\frac{\mathrm{rk}M} {\dim \bS_\lambda W}.$$ It follows from the proof of Theorem \ref{PieriBWBTheorem} that the rank of $M_k$ is the Euler characteristic of $H_\bullet(k)$, which is a polynomial in $k$ of degree $n-1$ with leading coefficient $\frac{\mathrm{rk} M}{(n-1)!},$ and that the dimension of $\bS_{\lambda^k}V$ as a function of $k$ is also a polynomial of degree $n-1$ with leading coefficient $\frac{\dim \bS_\lambda W}{(n-1)!} $. Hence, we have $$ \lim_{k\to \infty}\frac{\mathrm{rk} M_k}{\dim \bS_{\lambda^k }V}=\frac{\mathrm{rk}M }{\dim \bS_\lambda W}$$ as required. Thus, $\lim_{k\to \infty} C_k=C$ where $C_k$ is the corresponding lower bound obtained in the $\lambda^k, \mu^k$ case, and we are done.     
\end{proof}

\subsection{Generalizing to $\Sym^d V$ and $\bigwedge^d V$ when $d>2$}\label{generalisationSection}
In this section, we discuss some ideas on how we can construct partitions $\lambda$ and $\mu$ with non-maximal generic rank when $U=\Sym^d V$ or $\bigwedge^d V$ for $d>2$. 

The key step in the $\Sym^2 V$ and $\bigwedge^2 V$ cases was using a suitable Lie algebra and parabolic subalgebra whose nilpotent radical is given by $\Sym^2 V$ and $\bigwedge^2 V$ respectively. This allowed us to construct resolutions over the symmetric algebra using Kostant's theorem. Although it is not possible to find a suitable Lie algebra and parabolic with nilpotent radical equal to $U$ for general $U$, it is possible to construct a Kac-Moody algebra with a parabolic subalgebra such that its Levi subalgebra is given by $\fgl(V)$ and $U$ appears as the degree-one-graded piece of the nilpotent radical (as a representation of the Levi). 

Kostant's theorem allows us to construct free resolutions over the universal enveloping algebra of the whole nilpotent radical. Base-changing from the universal enveloping algebra of the nilpotent radical to $\Sym(U)$ (which is a quotient) gives us a complex over $\Sym (U)$ which is no longer exact. However, if the homologies of this complex are not supported at the generic point, the techniques in the proofs of Theorems \ref{Sym2Pieri} and \ref{Wedge2Pieri} still apply. This happens for example when the Kac-Moody algebra is a finite dimensional simple Lie algebra. Let us do an example to illustrate this point.

\begin{example}\label{E6 example}
    Let $\fg=E_6$, the exceptional simple finite dimensional Lie algebra with the following Dynkin diagram. 

    \begin{figure}[h]
    \centering
    \includegraphics[width=0.3\linewidth]{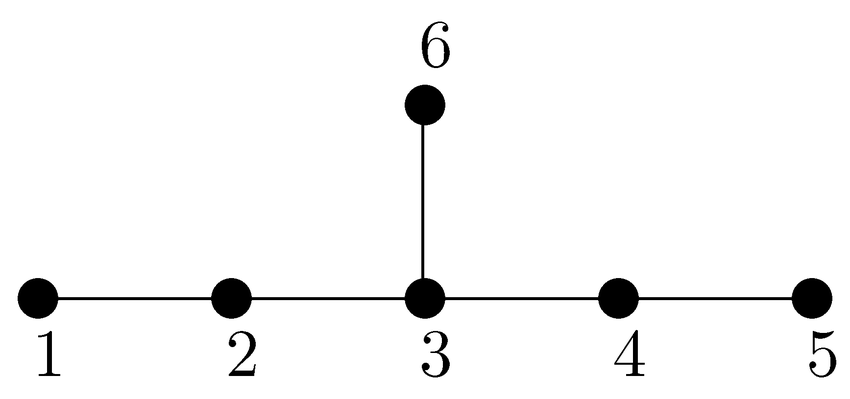}
    \label{fig:E_6}
\end{figure}

    We consider the parabolic subalgebra $\fp$ corresponding to the $6$th node according to the Dynkin diagram. We have a $\mathbb{Z}$-grading on the parabolic subalgebra 
    $$\fp=\fp_{-2}\oplus \fp_{-1}\oplus \fp_0\oplus\fp_1\oplus \fp_2.$$
    $\fp_0$ is the Levi subalgebra, which is isomorphic to $\fgl(\mathbb{C}^6)$. The nilpotent radical $\fu$ consists of the positive degree pieces. As a representation of $\fgl(\mathbb{C}^6)$, we have $\fp_1\cong \bigwedge^3\mathbb{C}^6$, which will be our $U$ for this example (and $V$ will denote $\mathbb{C}^6$), and $\fp_2\cong \mathbb{C}$, the trivial representation. 

    Since it is not very easy to describe the dominant weights and the Weyl group (and its action) of $E_6$, we will need to work abstractly. For $1\leq i\leq 6$, let $\alpha_i$ denote the simple roots of $E_6$, $h_i$ the corresponding simple coroots and $w_i$ the corresponding fundamental weights. The Weyl group $W$ of $E_6$ is generated by the simple reflections $s_i$ for $1\leq i\leq 6$. We recall that the action of the simple reflections on a weight $\beta$ is given by the formula $$s_i(\beta)=\beta-\langle \beta,h_i\rangle\alpha_i.$$
    The value of $\langle \alpha_i,h_j\rangle$ is given by the $(i,j)$th entry of the Cartan matrix: 
    
    $$\begin{bmatrix}
       2 & -1 & 0 & 0 & 0 & 0 \\ 
       -1 & 2 & -1 & 0 & 0 & 0 \\ 
       0 & -1 & 2 & -1 & 0 & -1 \\ 
       0 & 0 & -1 & 2 & -1 & 0 \\ 
       0 & 0 & 0 & -1 & 2 & 0 \\ 
       0 & 0 & -1 & 0 & 0 & 2  
    \end{bmatrix}$$

We also note that it follows from the definition of the $w_i$'s that $$s_i(w_j)=w_j-\delta_{ij}\alpha_i.$$

    Any integral dominant weight $\beta$ can be written as $$\beta=\sum_{i=1}^6a_iw_i$$ for a tuple $a=(a_1,\dots,a_6)$ of non-negative integers. We denote the corresponding finite dimensional irreducible representation by $L(\beta)$. We have $$\rho=\sum_{i=1}^6 w_i.$$ A fact that will be useful to us is that $$s_i(\rho)=\rho-\alpha_i.$$

    We are almost ready to apply Kostant's theorem for our setup. Let $W_6$ denote the Weyl group of the Levi and $W_6'$ the set of minimal length representatives of the cosets of $W/W_6$. We list the elements of $W_6'$ with lengths less than or equal to $3$ and the dotted action of their inverses on the dominant weight $\beta$ corresponding to the tuple $(a_1,\dots,a_6)$:

\vspace{0.2in}

    \begin{center}
    \begin{tabular}{|c|c|c|}
    \hline
       $l(w)$  & $w$ & $w^{-1}\bullet \beta$\\
         \hline
         $0$ & $\mathrm{id}$ & $(a_1,a_2,a_3,a_4,a_5,a_6)$\\
         $1$ & $s_6$ & $(a_1,a_2,1+a_3+a_6,a_4,a_5,-2-a_6)$\\
         $2$ & $s_3s_6$& $(a_1,1+a_2+a_3,a_6,1+a_3+a_4,a_5,-3-a_3-a_6)$\\
         $3$ & $s_4s_3s_6$& $(a_1,2+a_2+a_3+a_4,a_6,a_3,1+a_4+a_5,-4-a_3-a_4-a_6)$\\   
         $3$ & $s_2s_3s_6$& $(1+a_1+a_2,a_3,a_6,2+a_2+a_3+a_4,a_5,-4-a_2-a_3-a_6)$\\
         \hline
    \end{tabular}

\end{center}

\vspace{0.2in}
We stress the fact that in the table above, by $(\nu_1,\dots,\nu_6)$, we mean the weight corresponding to $\sum \nu_iw_i$. This need not be a weakly decreasing sequence, but we get the corresponding Schur functor $\bS_\gamma V$ for  $\gamma=(\gamma_1,\dots,\gamma_6)$ (up to a twist by the determinant) using the formula $$\gamma_i=\sum_{j=i}^6 \nu_j.$$ 
We show the computation for the length $1$ element (the computation for the rest can be done similarly).

\begin{align*}
    s_6\bullet\beta=s_6(\beta+\rho)-\rho&=s_6(\sum a_iw_i)+s_6(\rho)-\rho\\&=\beta-a_6\alpha_6+\rho-\alpha_6+\rho\\&=\beta-(a_6+1)\alpha_6
\end{align*}
to write $\beta-(a_6+1)\alpha_6$ in the form $\sum\nu_iw_i$, we use the fact that the coefficient of $w_i$ is $$\langle \beta-(a_6+1)\alpha_6,h_i\rangle.$$
Thus, when $i\neq 3,6$, using the Cartan matrix we observe that the coefficient of $w_i$ is $a_i$, and the coefficient of $w_3$ is $a_3+a_6+1$ while the coefficient of $w_6$ is $a_6-2(a_6+1)=-2-a_6$.

Thus, the first few terms of the resolution $\mathbb{F}_\bullet$ of $L(\beta)^*$ over the enveloping algebra $\fU$ of the nilpotent radical (up to a twist by the determinant) is:

$$\dots \to (\bS_{\beta_4}V\oplus \bS_{\beta_3}V)\otimes \fU \to \bS_{\beta_2}V\otimes \fU \to  \bS_{\beta_1}V\otimes \fU \to (\det V)^k\otimes \bS_{\beta}V^*\otimes \fU$$ 
    where 
       $$\beta_1=(b_1-b_6+a_6+1,b_1-b_5+a_6+1,b_1-b_4+a_6+1,b_1-b_3,b_1-b_2),$$ 
        $$\beta_2=(b_1-b_6+a_3+a_6+2,b_1-b_5+a_3+a_6+2,b_1-b_4+a_6+1,b_1-b_3+a_3+1,b_1-b_2),$$
        $$\beta_3=(b_1-b_6+a_3+a_4+a_6+3,b_1-b_5+a_3+a_6+2,b_1-b_4+a_4+a_6+2,b_1-b_3+a_3+a_4+2,b_1-b_2),$$
       $$\beta_4=(b_1-b_6+a_2+a_3+a_6+3,b_1-b_5+a_2+a_3+a_6+3,b_1-b_4+a_6+1,b_1-b_3+a_3+1,b_1-b_2+a_2+1).$$ 
Here, we use $b_i$ to denote the sums $$b_i=\sum_{j=i}^6 a_j.$$

We tensor $\mathbb{F}_\bullet$ with $\Sym U \cong \fU/\fU_{\fp_2}$ to obtain a complex over $\Sym U$. By definition, the homology of this complex is given by the Tor groups $\mathrm{Tor}_{\fU}(L(\beta)^*,\Sym U)$. Since $L(\beta)^*$ is a finite-dimensional vector space, it has proper support and hence the Tor groups have proper support as well. Consequently, the rest of the proof of Theorem \ref{Sym2Pieri} works just fine here, and we obtain families of pairs of partitions over $\bigwedge^3\mathbb{C}^6$ that are not generically surjective.

We observe from our description of the $\beta_i$'s that we unfortunately are not able to produce families where the boxes are added along a continuous border strip. However, we are able to produce several other families that are not generically surjective:
\begin{enumerate}
    \item If $a_3=0$, then $\beta_2$ is obtained from $\beta_1$ by adding a box to the $1$st, $2$nd and $4$th rows.
    \item If $a_4=0$, then $\beta_3$ is obtained from $\beta_2$ by adding a box to the $1$st, $3$rd and $4$th rows.
    \item If $a_2=0$, then $\beta_4$ is obtained from $\beta_2$ by adding a box to the $1$st, $2$nd and $5$th rows.
\end{enumerate}
If one could show that the Young flattening for these pair of partitions is of full rank (or sufficiently large), then that would give non-trivial lower bounds for tensors coming from these Pieri maps.
\end{example}

In general, when $\fg$ is not a finite dimensional simple Lie algebra, the irreducible modules $L(\beta)$ are no longer finite dimensional. Thus, it is not clear if the Tor modules $\mathrm{Tor}(L(\beta),U)$ have finite support, which is an obstacle for generalizing the methods of this paper for general $U$. Further, as we see in Example \ref{E6 example}, we may not be able construct exact sequences using Kostant's Theorem that contain the Pieri map (as a differential) between $\lambda$ and $\mu$ where the partitions differ by adding boxes to the same row or column, which is the case for which we know that the Young flattening is of full rank.